\documentclass[twocolumn]{article}

\usepackage[affil-it]{authblk}

\usepackage{rotating}

\usepackage[short,nocomma]{optidef}

\usepackage[latin1]{inputenc} 

\usepackage{mathtools}

\usepackage[english]{babel}

\usepackage{booktabs}

\usepackage{pdflscape}

\usepackage{enumitem}

\usepackage[framemethod=tikz]{mdframed}
\newmdenv[innerlinewidth=0.5pt, roundcorner=4pt,linecolor=black,innerleftmargin=6pt,
innerrightmargin=6pt,innertopmargin=6pt,innerbottommargin=6pt]{mybox}

\usepackage{ mathrsfs }

\usepackage[export]{adjustbox}

\usepackage{chngpage}

\usepackage[authoryear,longnamesfirst]{natbib}

\usepackage{multirow}

\usepackage{amssymb}
\usepackage{amsmath}

\usepackage{bbm}

\usepackage{array}
\usepackage{booktabs}
\usepackage{colortbl}

\usepackage{amsthm}

\newtheorem*{myprop*}{Proposition}

\usepackage{float}

\usepackage[hyphens]{url}

\usepackage{fullpage}

\usepackage[colorinlistoftodos]{todonotes}

\usepackage{xcolor,soul}
\definecolor{dr}{rgb}{0.75,0.00,0.00}
\definecolor{lr}{rgb}{1.00,0.75,0.75}
\sethlcolor{lr}

\definecolor{dr2}{rgb}{0.00,0.00,0.00}
\definecolor{lr2}{rgb}{0.00,1.00,0.00}

\definecolor{dr3}{rgb}{0.00,0.00,1.00}
\newcommand{\changed}[1]{#1}
\usepackage{algorithm}
\usepackage{algpseudocode}
\algdef{SE}[DOWHILE]{Do}{doWhile}{\algorithmicdo}[1]{\algorithmicwhile\ #1}

\usepackage{floatpag}
\floatpagestyle{plain}

\usepackage{subfig}

\usepackage{afterpage}

\usepackage{tikz}
\usetikzlibrary{decorations.pathmorphing, arrows.meta, shapes, backgrounds, decorations.markings}
\usetikzlibrary{arrows, automata}

\usepackage{graphicx}

\usepackage{longtable}

\usepackage{threeparttable}

\usepackage{arydshln}

\usepackage{pdflscape}

\newenvironment{nospacebelowflalign*}
 {\setlength{\belowdisplayskip}{0pt}%
  \csname flalign*\endcsname}
 {\csname endflalign*\endcsname\ignorespacesafterend}

\newenvironment{nospaceflalign*}
 {\setlength{\abovedisplayskip}{0pt}\setlength{\belowdisplayskip}{0pt}%
  \csname flalign*\endcsname}
 {\csname endflalign*\endcsname\ignorespacesafterend}

\newcolumntype{L}[1]{>{\raggedright\let\newline\\\arraybackslash\hspace{0pt}}m{#1}}
\newcolumntype{C}[1]{>{\centering\let\newline\\\arraybackslash\hspace{0pt}}m{#1}}
\newcolumntype{R}[1]{>{\raggedleft\let\newline\\\arraybackslash\hspace{0pt}}m{#1}}

\usepackage{caption} 
\graphicspath{ {images/} }

\usepackage[colorlinks = true,
            linkcolor = blue,
            urlcolor  = blue,
            citecolor = black,
            anchorcolor = blue]{hyperref}

\title{Resiliency of On-Demand Multimodal Transit Systems \\ During a Pandemic}

\author[1,2]{Ramon Auad}
\author[1]{Kevin Dalmeijer}
\author[1]{Connor Riley}
\author[1]{Tejas Santanam}
\author[1]{Anthony Trasatti}
\author[1]{Pascal Van Hentenryck}
\author[1]{Hanyu Zhang}
\affil[1]{\normalsize{H. Milton Stewart School of Industrial and Systems Engineering, Georgia Institute of Technology}}
\affil[2]{\normalsize{Departamento de Ingenier\'ia Industrial, Universidad Cat\'olica del Norte}}

\date{\today}

\begin{document}






\pagestyle{plain}
\pagenumbering{arabic}


\twocolumn[
  \begin{@twocolumnfalse}
    \maketitle

\begin{abstract}
\footnotesize
\noindent
During the COVID-19 pandemic, the collapse of the public transit
ridership led to significant budget deficits due to dramatic decreases
in fare revenues.  Additionally, public transit agencies are facing
challenges of reduced vehicle capacity due to social distancing
requirements, additional costs of cleaning and protective equipment,
and increased downtime for vehicle cleaning.  Due to these constraints
on resources and budgets, many transit agencies have adopted essential
service plans with reduced service hours, number of routes, or
frequencies. 
This paper studies the resiliency during a pandemic of On-Demand
Multimodal Transit Systems (ODMTS), a new generation of transit
systems that combine a network of high-frequency trains and buses with
on-demand shuttles to serve the first and last miles and act as
feeders to the fixed network.  It presents a case study for the city
of Atlanta and evaluates ODMTS for multiple scenarios of depressed
demand and social distancing representing various stages of the
pandemic.

The case study relies on an optimization pipeline that provides an end-to-end ODMTS solution by bringing together methods for demand estimation, network design, fleet sizing, and real-time dispatching.
These methods are adapted to work in a multimodal setting and to satisfy practical constraints. In particular, a limit is imposed on the number of passenger transfers, and a new network design model is introduced to avoid the computational burden stemming from this constraint. Real data from the Metropolitan Atlanta Rapid Transit Authority (MARTA) is used to conduct the case study, and the results are evaluated with a high-fidelity simulation. The case study demonstrates how ODMTS provide a resilient solution in terms of cost, convenience, and accessibility for this wide range of scenarios.

\end{abstract}

  \end{@twocolumnfalse}
]


\section{Introduction}
\label{sec:introduction}
At the start of the COVID-19 pandemic, public transit ridership collapsed across the nation as many individuals stayed at home. 
Public transit agencies saw their budget deficits build as ridership dropped as much as 85\% in some cities, leading to a dramatic drop in fare revenue for 2020 compared to 2019 \citep{GoldbaumWright2020-existentialPerilMass}. In addition, public transit agencies are facing the challenging task of providing safe and accessible travel, especially for essential workers, with reduced vehicle capacities, additional costs for protective equipment or supplies, and increased downtime for cleaning. 

With constrained vehicle capacity and large budget deficits, many transit agencies have adopted essential service plans that include cutting routes, service hours, or frequencies in order to provide safer service on other routes with their limited resources.
New York City's MTA no longer runs 24 hours a day (it closes from 1-5am to provide downtime for more effective cleaning) and faces an \$8B deficit \citep{Nessen2020-NewYorkCitys}. In Boston, the MBTA has tried to preserve as much of the weekday trips as possible for essential workers, but is now looking at dramatic budget cuts due to a deficit of close to \$600M.  \citep{GoldbaumWright2020-existentialPerilMass}. At the time of writing, the MBTA is considering cutting their ferry routes and weekend commuter rail services. In Atlanta, the Metropolitan Atlanta Rapid Transit Authority (MARTA) faced similar challenges. In April, MARTA reduced the frequency of the rail network from every 12 minutes to every 20 minutes and its essential service plan reduced the number of bus lines from 110 to 41 \citep{Wickert2020-MartaCutsMost}.
Even as transit planners and operators try to mitigate the effects of these service cuts, many essential workers and community members are still impacted in terms of accessibility and convenience.
In Atlanta, many individuals rely on public transit for access to jobs, groceries, education, and healthcare, and are now forced to find alternative forms of transit \citep{Wickert2020-MartaBusCuts}.
Ridesharing services such as Uber and Lyft are often prohibitively expensive, which makes it crucial that transit services return to normal as soon as possible.

\begin{figure}[!t]
    \centering
    \includegraphics[width=\linewidth]{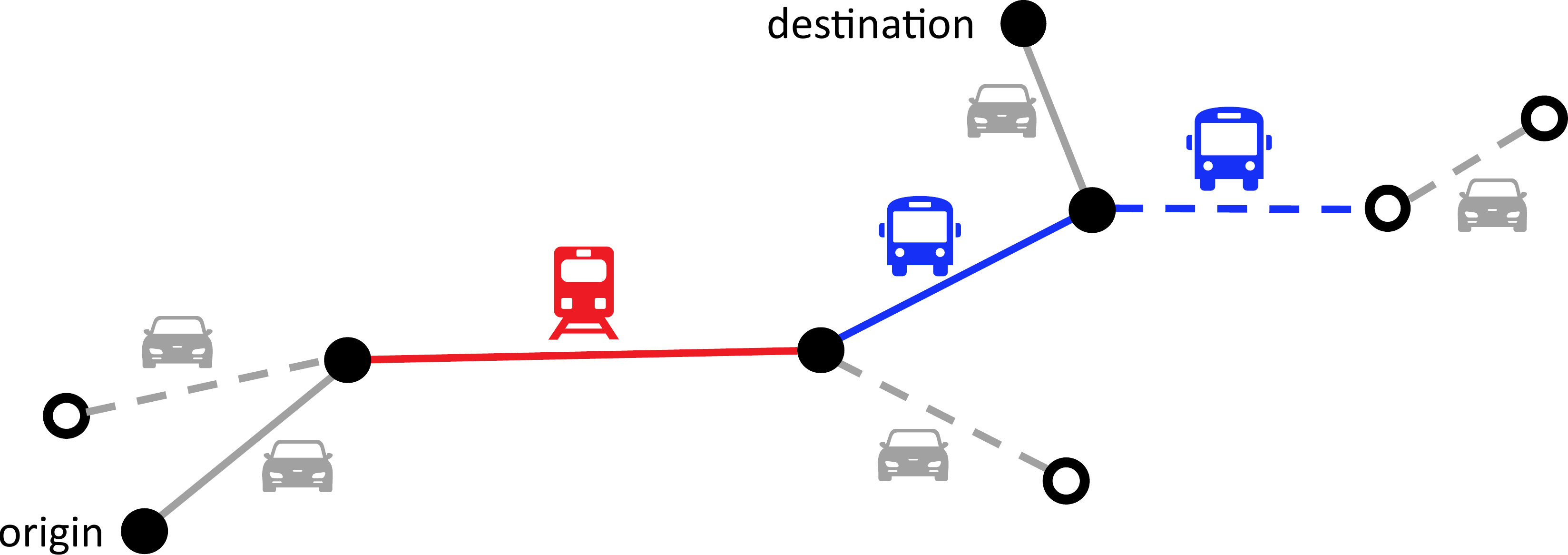}
    \caption{An ODMTS and an Example \changed{of a} Passenger Path (solid lines).}
    \label{fig:ODMTS-Example}
\end{figure}

This paper studies the resiliency of On-Demand Multimodal Transit Systems (ODMTS) \citep{MaheoEtAl2019-BendersDecompositionDesign} during a pandemic.
ODMTS is a new generation of transit systems that consist of a network of
high-frequency trains and/or buses combined with on-demand shuttles to
serve the first and the last miles. In this paper, the \emph{resiliency} of an ODMTS is defined as the measure in which the system is able to provide accessible public transit with a high level of service, within an appropriate budget, when a pandemic results in significantly lower ridership and additional safety requirements. ODMTS can respond to such a decrease in ridership by scaling down the number of on-demand shuttles and/or high-frequency buses.
On-demand shuttles keep serving passengers close to where they are, such that access to the system is not compromised.
By design, trains and buses are only used to serve busy
corridors, which would still see traffic during a
pandemic. Figure~\ref{fig:ODMTS-Example} illustrates the design of an
ODMTS and how it operates. This ODMTS features fixed rail and rapid
bus transit routes, as well as on-demand shuttles. While the trains
and buses operate on a fixed schedule, the shuttles are routed
dynamically in real time. The example path through the ODMTS involves
a shuttle leg to pick up the rider, a rail leg, a bus leg, and a final
shuttle leg to bring the rider to her destination. All these legs are
synchronized: for instance, the shuttle serving the last leg is
waiting at the bus stop to bring the rider to her final
destination. This
\href{https://sam.isye.gatech.edu/projects/demand-multimodal-transit-systems}{video}
describes this process in more detail \citep{SAML2020-HowDemandMultimodal}.

Real-world case studies of ODMTS systems have repeatedly found that ODMTS have the potential to cut costs and improve passenger convenience.
\citet{MaheoEtAl2019-BendersDecompositionDesign} present a case study for Canberra, Australia, for which the current network comprises about 2,800 bus stops, and which sees about 60,000 passengers per day.
ODMTS simulations have been performed for the transit system of the University of Michigan, and some of the results have been validated by a pilot in the Spring of 2018 \citep{VanHentenryck2019-DemandMobilitySystems}.
\citet{BasciftciVanHentenryck2020-BilevelOptimizationDemand, BasciftciVanHentenryck2021-CapturingTravelMode} and \citet{AuadVanHentenryck2021-RidesharingFleetSizing} study the Ann-Arbor/Ypsilanti region in Michigan, which currently consists of 1,267 bus stops, and the case studies consider about 3,000 up to 7,000 passengers at a time.
\citet{DalmeijerVanHentenryck2020-TransferExpandedGraphs} present results for 7,167 riders of the MARTA system in Atlanta, which currently consists of 5,563 bus stops and 38 rail stations.
Aside from cutting cost and improving passenger convenience, other advantages of ODMTS and its contribution to social mobility are discussed by \citet{VanHentenryck2019-DemandMobilitySystems}.

To study the resiliency of ODMTS to the effects of depressed ridership and additional safety requirements that arise during a pandemic, an ODMTS is designed for the city of Atlanta, and its performance is evaluated under different scenarios that capture various stages of the pandemic.  The evaluation is based on a detailed simulation that uses real transit data provided by MARTA. \changed{The contributions of the paper can be summarized as follows:}
\begin{enumerate}
	\item \changed{The paper presents a pipeline for designing and evaluating an ODMTS, combining advanced models for trip chaining, ODMTS design, ridesharing and fleet sizing, and real-time shuttle dispatching for the first time. These models are adapted to operate in a multimodal setting, and to satisfy practical constraints. In particular, the models impose a constraint on the number of passenger transfers, and a new network design model is introduced to avoid the computational burden imposed by this constraint.}
	\item \changed{The paper provides a comprehensive case study for the city of Atlanta to assess the resiliency of ODMTS under various scenarios of pandemic response, which include loss of revenue due to depressed demand and reduced vehicle capacity due to social distancing requirements. High-fidelity simulations show that, by scaling down the number of shuttles appropriately, an ODMTS provides a similar level of service with a cost savings that compensates for revenue losses. The case study also shows that an ODMTS can operate successfully under more restrictive budgets, at the cost of increasing waiting times during the busiest times of the day.}
\end{enumerate}

This rest of the paper is structured as follows.  Section~\ref{sec:related_work}
presents related work in the literature, and
Section~\ref{sec:reviewmarta} reviews the current MARTA transit system
and its response to the pandemic.
The pipeline for designing and evaluating an ODMTS is presented in Section~\ref{sec:method}.
Section~\ref{sec:baseline} introduces an ODMTS for Atlanta, and its resiliency is tested in Section~\ref{sec:casestudy}.  The final
section summarizes the findings and presents directions for future
research.

\section{Related Work}
\label{sec:related_work}

Recent work has focused on understanding the impact of COVID-19 on
public transit ridership.  \citet{LiuEtAl2020-ImpactsCovid19} provide a
systematic analysis of 113 public transit systems across the U.S.,
using data from the Transit App to study the unprecedented decline of
public transit during this pandemic.  They explored how social
distancing, self-quarantine, and working from home recommendations
from the Center for Disease Control and Prevention (CDC) impacted
ridership, even before the local spread of COVID-19.
\citet{TirachiniCats2020-Covid19Public} discuss factors that influence the contagion risk of COVID-19 in public transit, including vehicle and station occupancy, exposure time, mask use, and hygiene.
COVID-19 has led to new challenges for the transportation community, and different examples are given in the survey by \citet{AgatzEtAl2021-MakeNoLittle}.
\citet{GkiotsalitisCats2021-PublicTransportPlanning} study the impact of COVID-19 on public transit, and survey methods for public transit planning that can be used to address the changes in demand patterns and reductions in capacity.
These include methods to set service frequencies, modify timetables, and to prevent crowding.
Changing the service frequencies depending on the social distancing requirements has been studied by \citet{GkiotsalitisCats2021-OptimalFrequencySetting}.
The authors develop a mixed-integer quadratic programming model, which is used to conduct a case study for the Washington D.C. metro network.

The interaction between multimodal networks and
on-demand vehicles has been considered by
\citet{SalazarEtAl2018-InteractionAutonomousMobility}, who look into the coordination of
public transit systems with fleets of autonomous vehicles and its
benefits. Their case study is based on New York City and the existing
public transit system.  \citet{PintoEtAl2020-JointDesignMultimodal} focus on the allocation
of resources between public transit fleets and shared-use autonomous
mobility services.  They choose the frequency of the public transit
routes and size of the fleet in a bilevel problem under a constrained
budget.  A case study is provided for the city of Chicago.
\section{The MARTA Transit System}
\label{sec:reviewmarta}
%

\begin{figure}[!t]
    \centering
    \includegraphics[width=0.4\textwidth, trim=222 360 175 220, clip]{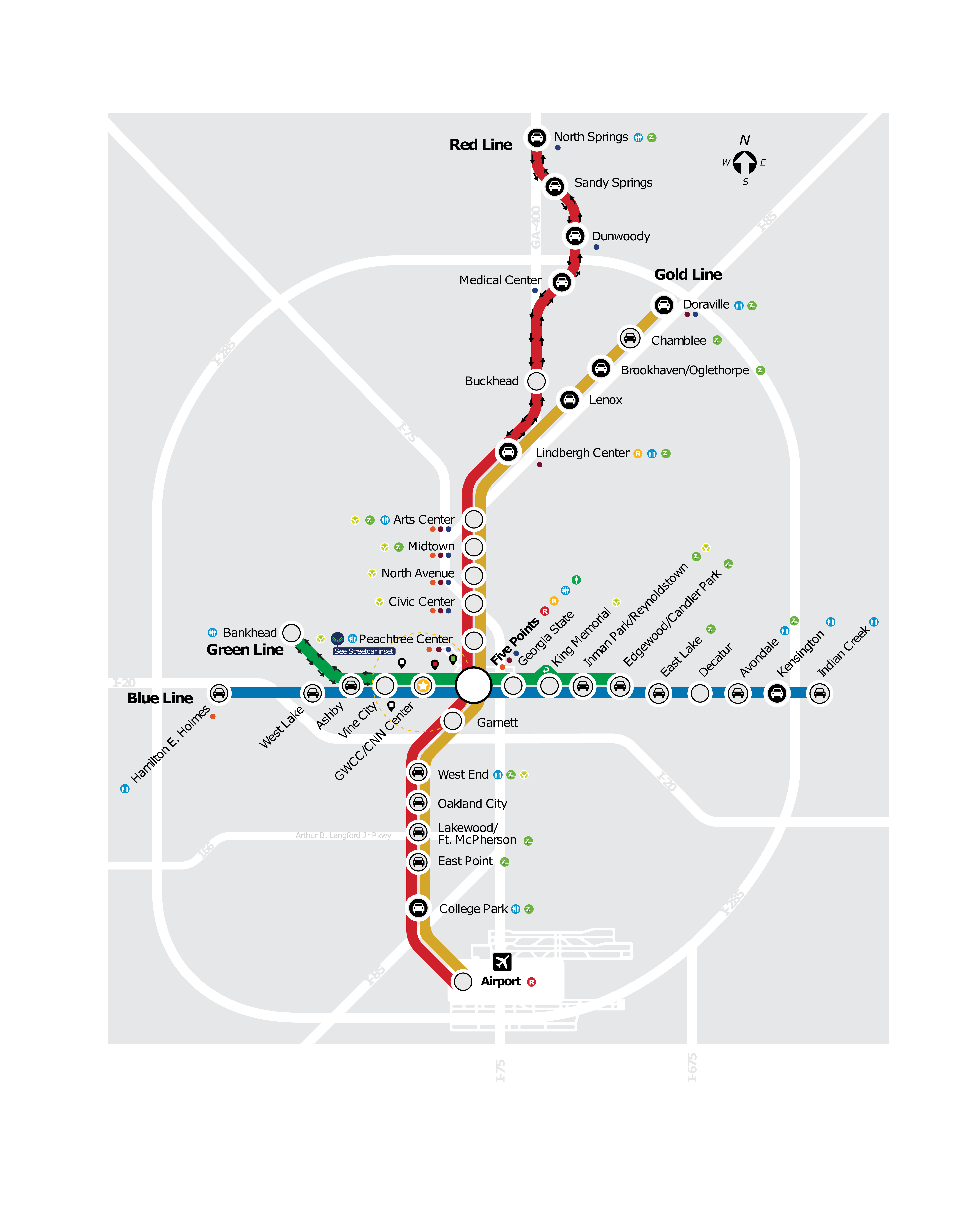}
    \caption{The MARTA Rail Network \citep{MARTA2020-TrainStations}.}
    \label{fig:marta_rail}
\end{figure}

\begin{figure*}[!t]
    \centering
    \begin{minipage}{0.5\linewidth}
        \centering
        \includegraphics[width=0.9\linewidth, trim=1180 0 1150 0, clip]{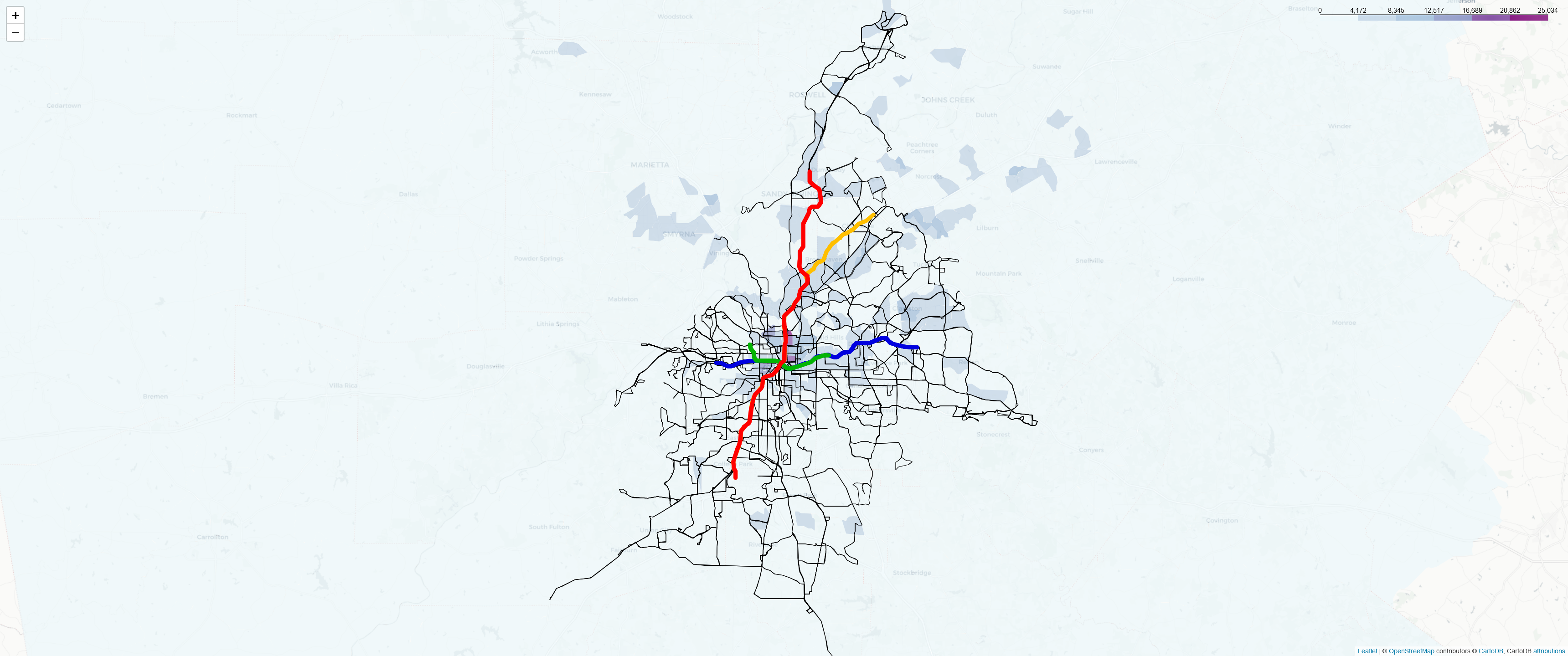}
        \caption{The MARTA Bus and Rail Coverage in 2019.}
        \label{fig:marta_network}
    \end{minipage}%
    \hfill
    \begin{minipage}{0.5\linewidth}
        \centering
        \includegraphics[width=0.9\linewidth, trim=1180 0 1150 0, clip]{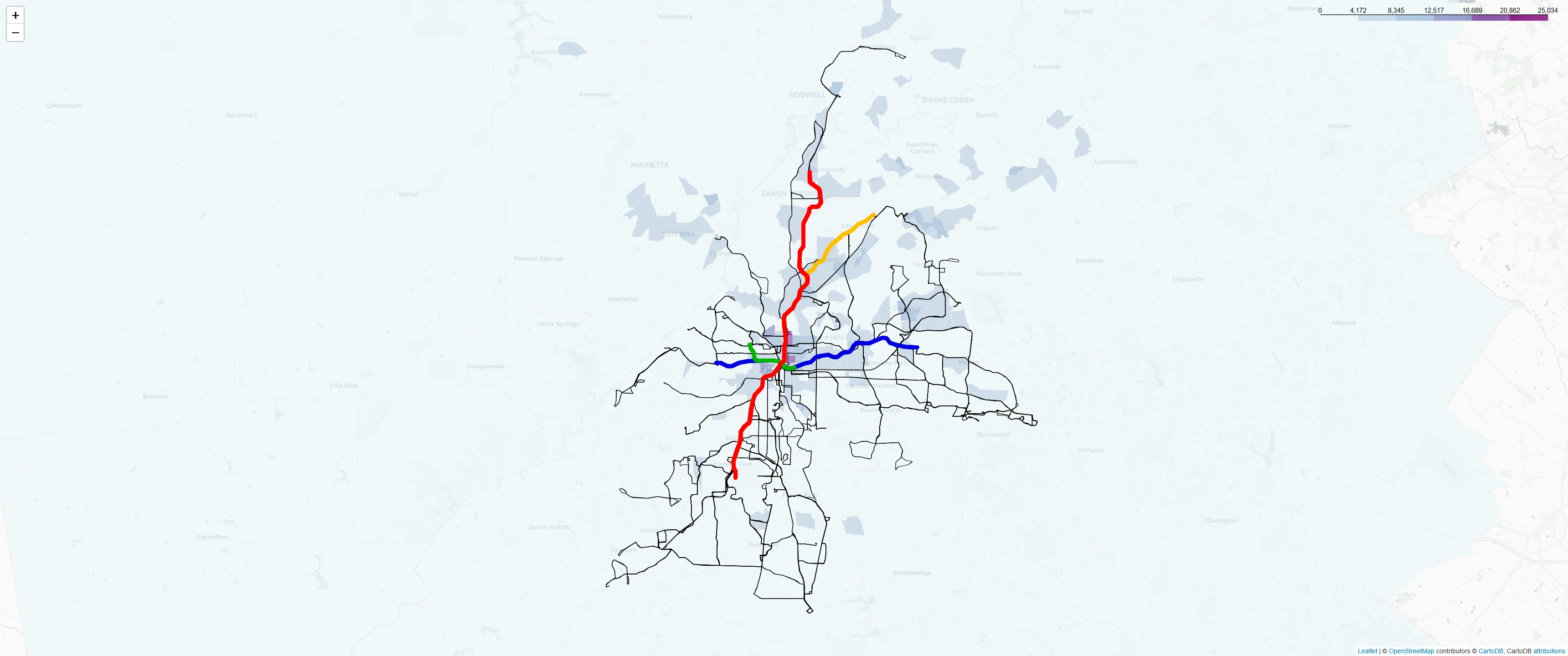}
        \caption{The MARTA Essential Service Plan.}
        \label{fig:marta_esp}
    \end{minipage}
\end{figure*}

The MARTA transit system consists of trains and buses that cover
Fulton, Clayton, and DeKalb counties in the metro Atlanta area.  The
rail network consists of 38 stations and four lines: the Green, Blue,
Red, and Gold line (Figure~\ref{fig:marta_rail}).  With around 160,000
transactions on a typical weekday, this is the 8th largest rail
network in the US by ridership.  In addition to the rail, MARTA
operates over 500 buses that serve over 110 different lines, and see
approximately 110,000 transactions every day.
Figure~\ref{fig:marta_network} shows the total bus and rail coverage,
the rail being the backbone of the system.
The network is overlayed on a map of population densities \citep{ARC2019-PopulationbyCensus}.

According to the \citet{CensusBureau2020-MeansTransportationWork},
the rate of public transit commuting to work in Atlanta is 10.4\%; in
contrast, NYC and Boston have rates of 56.0\% and 33.3\%,
respectively. Atlanta is also routinely ranked among the most
congested cities in the world.  The INRIX Traffic
Scorecard ranks Atlanta to be the 10th most congested city in the US,
with 82 lost hours per year for the average driver
\citep{Reed2019-GlobalTrafficScorecard}.  Still, it is reported that
driving typically results in a shorter commute than using public
transit.  This \changed{is in line} with Figure~\ref{fig:marta_network}, which shows
that the public transit network spans a large geographical area, but is also sparse in many regions, which prevents convenient access.

\subsection{Ridership Data}
\label{sec:reviewmarta:data}

To access the MARTA system, riders use a transit pass, known as the
\emph{Breeze Card}.  Passengers tap their Breeze Card when boarding a
bus and when entering and exiting a rail station. Inside buses,
riders can also pay with cash. The Breeze Card \emph{transaction data} is
collected by the Automated Fare Collection (AFC) system, which records
the time and type of the transaction, the terminal at which the
transaction took place, and the Breeze card id (if applicable).  In
addition to transaction data, MARTA collects \emph{Automated Passenger
  Counter (APC) data} in the MARTA buses, using cameras or infrared
beams in the doorways to record boardings and alightings.  For every
stop, the data contains the number of boardings and alightings, the
vehicle id, the route and direction, and the time and duration of the
stop.  This data is supplemented with GPS coordinates from the
Automatic Vehicle Location system.

\begin{figure}[!t]
	\centering
	\includegraphics[width=\linewidth, trim=0cm 0.1cm 3cm 1cm, clip]{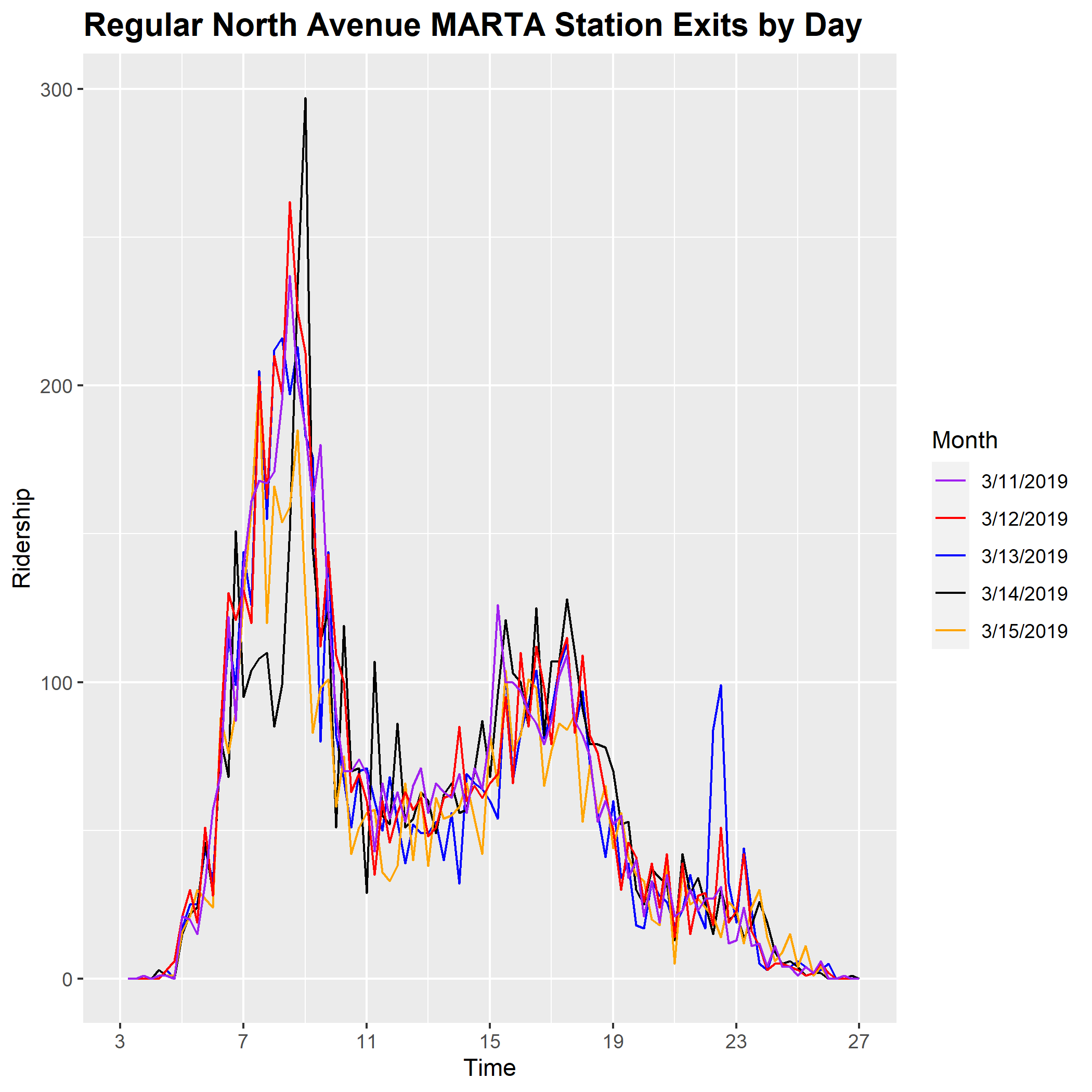}
	\caption{Arrivals at North Avenue Station per 15 Minute Interval (first five weekdays of March 2019).}
	\label{fig:nave_signature}
\end{figure}

For the case study, transaction data and APC data for 2018 are
combined to create a set of 192,311 passenger trips that is
representative for an average weekday in March 2018.  Details on
estimating the origin-destination matrix are provided in
Section~\ref{sec:method:odmatrix}.  The case study focuses on the
55,871 passengers traveling during the 6am--10am morning peak, which
is one of the most challenging periods of the day.  Transaction data
collected during the pandemic has been made available for March up to
May 2020.  To study the effect of the pandemic on ridership, the same
data was provided for the year before.
Ridership over the day is typically very predictable.
As an example, Figure~\ref{fig:nave_signature} presents the number of passengers arriving at North Avenue station per 15 minute interval for consecutive weekdays in March 2019, which shows that demand under normal circumstances is quite stable, even for short time intervals and at the station level.

\subsection{Cost and Performance}
\label{sec:costcurrent}

The \citet{FTA2018-NationalTransitDatabase} provides details on
MARTA's financials.  In 2018, operating and capital expenses were
\$591M and \$244M, respectively.  Most operating costs stem from the
bus network (49\%) and the rail (44\%); the remainder is for
paratransit services for passengers. Fares and directly generated
funds, which are most affected during a pandemic, pay for 33\% of the
operating budget.

MARTA reports an average operating expense of \$104.10 per bus revenue
hour, which includes costs for general administration and facility
maintenance, but excludes depreciation of the vehicles.  For the case
study, it is important to know how the cost of operating a bus
compares to the cost of operating a shuttle.  To present a fair
comparison, this paper focuses only on the main cost differentiators:
labor cost, vehicle maintenance cost, and depreciation cost.  This
puts the hourly cost of a bus at \$72.15.  MARTA operates 465 buses in
maximum service, which results in a daily cost of \$134k during the
four hour morning peak.  Details on this cost calculation are provided
by Appendix~\ref{app:costs}.

In terms of performance, the data shows that waiting times are
relatively long, especially for bus users.  The average duration for
morning-peak trips that include a bus is 46 minutes: 20 minutes of
travel time and 26 minutes of waiting time.  The travel time is
calculated from the transaction data, and the waiting time is
estimated based on the service frequencies.  Even without traffic, the
average trip duration would still be 40 minutes for these passengers.
The situation is better for riders only using rail: their trips take
25 minutes on average. Note that the time to get to the bus stop or
rail station is not included in this calculation, which means that the
actual trips take longer.

\subsection{Pandemic Response}
\label{sec:demand}
In March 2020, COVID-19 cases were increasing rapidly.  MARTA responded
by implementing a series of hygiene protocols and by reducing the
maximum occupancy of buses and trains, following new guidelines by the
CDC.  The following is a brief timeline of events during March and
April 2020:
\begin{itemize}
    \item March 11th: the World Health Organization characterized the novel coronavirus outbreak as a pandemic \citep{WHO2020-WhoDirectorGenerals};
    \item March 26th: MARTA required passengers to enter buses through the rear door, and fare collection was suspended to accommodate this;
    \item March 30th: MARTA reduced the frequency of its bus schedule;
    \item April 20th: MARTA switched to an \emph{essential service plan} consisting of 41 bus routes.
\end{itemize}
Additionally, the rail frequencies were reduced during this period \citep{ARC2020-HowCovid19}.

\begin{figure}[!t]
	\centering
	\includegraphics[width=0.45\textwidth]{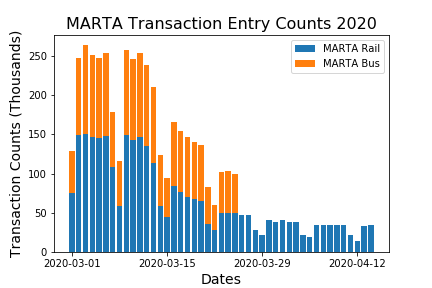}
	\caption{Breeze Card Tap-Ins from March to Mid-April 2020 (Bus transactions were not recorded after March 26).}
	\label{fig:RidershipComp:2020}
\end{figure}

Following the spread of COVID-19, public transit ridership collapsed
due to lockdowns and social distancing, and many people have started
working from home.  Figure~\ref{fig:RidershipComp:2020} shows how the
number of bus and rail transactions decreased during March and April
2020.  The ridership at the beginning of March was similar to that of
2019, and represented a typical demand.  Late March, bus ridership had
decreased by roughly 50\%, and rail ridership was reduced by
approximately 70\%.  On March 26, bus fare collection was suspended,
and no more transactions were recorded.  In April, rail ridership was
at around 22\% of the normal level, and bus ridership was down by 51\%
on April 19 (see
\href{https://www.itsmarta.com/KPIRidership.aspx}{MARTA performance
  indicators}). Ridership remained low thereafter, with similar
numbers for May.

MARTA's essential service plan was introduced to provide safer service on its most
critical routes, while dealing with rapidly declining ridership, additional cleaning 
requirements, and new social distancing policies.
It cut many routes that were assessed as ``less-critical''.
Figure~\ref{fig:marta_esp} shows the bus and rail coverage of the
essential service plan, and compared to
Figure~\ref{fig:marta_network}, it is clear that public transit is
less accessible during the pandemic.  Although this was necessary,
switching to the essential service plan did have a significant
negative impact on many people of Atlanta
\citep{Wickert2020-MartaBusCuts}.

\section{ODMTS}
\label{sec:method}

\begin{figure*}[!p]
    \centering
    \includegraphics[scale=0.65]{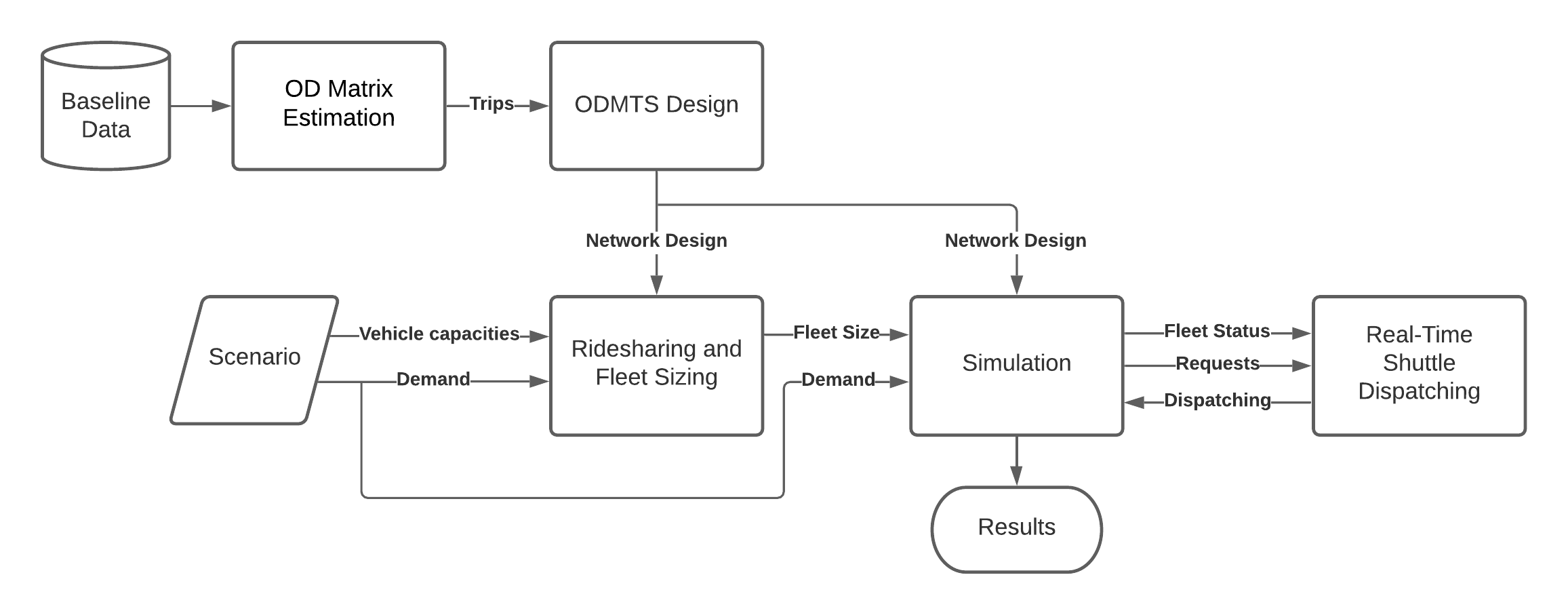}
    \caption{\changed{Pipeline for the Design and Evaluation of ODMTS.}}
    \label{fig:flowchart}
    \vspace{0.5 cm}
    \scriptsize
    \centering
    \begin{tabular}{p{0.1\textwidth}p{0.85\textwidth}}
        \toprule
        Symbol & Definition \\
        \midrule
        \textbf{Sets} & \\
        $V$ & set of locations\\
        $V_H$ & set of transit hub locations, $V_H \subseteq V$\\
        $M$ & set of modes, $M=\{shuttle, bus, rail\}$\\
        $F^m$ & set of frequencies for mode $m\in \{bus,rail\}$ (possible number of vehicles over the time horizon)\\
        $A$ & set of arcs, each $a \in A$ represents traveling from $i(a) \in V$ to $j(a) \in V$ with mode $m\in M$ and frequency $f(a) \in F^m$\\
        $A^{bus}$ & set of bus arcs, $A^{bus} \subseteq A$\\
        $P$ & set of shuttle requests, each $p \in P$ consists of an origin, destination, and request time\\
        $T$ & set of trips, each $t\in T$ consists of an origin $o(t) \in V$, destination $d(t) \in V$, number of passengers $p(t)$\\
        $P^+_h$ & shuttle requests with hub $h \in V_H$ as the origin, $P^+_h \subseteq P$\\
        $P^-_h$ & shuttle requests with hub $h \in V_H$ as the destination, $P^-_h \subseteq P$\\
        $P_{direct}$ & shuttle requests that do not involve a hub, i.e., $P_{direct} = P \setminus P_H$\\
        $R_h^+$ & set of shuttle routes that only serve requests in $P^+_h$\\
        $R_h^-$ & set of shuttle routes that only serve requests in $P^-_h$\\
        $S$ & set of shuttles, each $s \in S$ is associated with a current time and location\\
        $R_s$ & set of shuttle routes that can be performed by shuttle $s \in S$\\
        $R$ & shuttle routes considered by the dispatching algorithm, $R = \bigcup_{s\in S} R_s$\\
        \multicolumn{2}{l}{\textbf{Parameters}} \\
        $\alpha$ & weight parameter that assigns weight $\alpha$ to total trip duration and $(1-\alpha)$ to system cost, $\alpha \in [0,1]$\\
        $\tau_a$ & time to travel arc $a \in A$\\
        $d_a$ & distance to travel arc $a \in A$\\
        $c^{m}$ & hourly cost of operating one vehicle of mode $m\in M$\\
        $\bar{c}^{shuttle}$ & estimated shuttle cost per mile used by the network design algorithm\\
        $L$ & length of the time horizon\\
        $L^+$ & length of the extended time horizon used to get accurate results over the time horizon\\
        $\beta_a$ & cost of including arc $a \in A^{bus}$ in the design, $\beta_a = (1-\alpha) \tau_a f(a) c^{bus}$\\
        $\gamma_a^t$ & cost of traversing arc $a\in A$ for trip $t \in T$, $\gamma_a^t = p(t)\left((1-\alpha) d_a \bar{c}^{shuttle} + \alpha \tau_a\right)$ for shuttle, $\gamma_a^t = p(t)\alpha\left(\tau_a + \frac{L}{2f(a)}\right)$ for bus and rail\\
        $K$ & maximum number of arcs travelled by each passenger, $K \ge 1$, corresponds to $K-1$ transfers\\
        $\gamma^r$ & cost of shuttle route $r \in R_H$ as used by the ridesharing algorithm\\
        $\rho$ & tolerance factor for shuttle trips in the ridesharing algorithm\\
        $\Delta$ & time window width to allow for ridesharing in the ridesharing algorithm\\
        $b^r_p$ & indicator that is one if shuttle route $r$ covers shuttle request $p \in P$, and zero otherwise\\
        $Q^m$ & maximum number of passengers in one vehicle of mode $m\in M$\\
        $l$ & length of a dispatching epoch\\
        $c_r$ & cost of shuttle route $r \in R$ as used by the dispatcher, equal to sum of the waiting times\\
        $g_p$ & penalty for not serving request $p \in P$ in the current epoch\\
        $\bar{g}$ & initial penalty for not serving request $p\in P$, based on zero waiting time\\
        \multicolumn{2}{l}{\textbf{Variables}} \\
        $z_a$ & one if arc $a\in A^{bus}$ is part of the network design, zero otherwise\\
        $y_a^t$ & one if trip $t\in T$ uses arc $a \in A$ in the network design, zero otherwise\\
        $x_r$ & one if shuttle route $r$ is selected in the ridesharing solution, zero otherwise\\
        $f_{(r,r')}$ & one if shuttle routes $r, r' \in R_H$ are served sequentially by the same shuttle in the fleet-sizing algorithm, zero otherwise\\
        $\theta_r$ & one if shuttle route $r \in R$ is selected by the dispatcher, zero otherwise\\
        $w_p$ & one if shuttle request $p \in P$ is unserved in the current solution by the dispatcher, zero otherwise\\
        \bottomrule
    \end{tabular}
    \captionof{table}{Nomenclature.}%
    \label{tab:symbols}%
\end{figure*}

This section reviews the software pipeline for designing an ODMTS and
evaluating its cost and performance.
The pipeline provides an end-to-end solution by bringing together recent methods for demand estimation, network design, fleet sizing, real-time dispatching, and simulation, and adapting them to work in a multimodal setting.
Figure~\ref{fig:flowchart}
presents an overview of the pipeline, and Table~\ref{tab:symbols}
summarizes the nomenclature used in this section.

The {\em network design} of an ODMTS takes as input an Origin-Destination (O-D) matrix
which is estimated from the \emph{baseline data} and represents rider
trips, as well as a set of hubs that can be connected by rail or
buses. The design indicates which connections are served by buses,
rail, or shuttles, as well as the frequencies for the fixed lines.
The ridesharing and fleet-sizing optimizations then
determine the necessary number of shuttles. Based on the network
design, the fleet size, and the demand, a {\em real-time simulation}
mimics real-time operations. It keeps track of every rider and every
vehicle over time, including their positions over time, the trip legs
for riders, the vehicle occupancies, and the unserved requests.  The
real-time dispatching and routing of the shuttles is performed using
an optimization algorithm over a rolling horizon: the optimization is
run every 30 seconds on requests that have been batched in the
previous time periods or have not been served yet.
This leads to a high-fidelity evaluation of the ODMTS.

The simulator obviously depends on its inputs, but is independent from the models and assumptions that are used to generate these inputs.
In Section~\ref{sec:method:design}, for example, vehicle capacity are ignored while creating the design, but the simulator enforces that passengers wait for the first vehicle with an empty seat. Modeling choices are thus made during network design and fleet sizing, but they do not impact how the system is evaluated.
Capacity could be incorporated in the network design, but this is not needed in practice, since the large vehicles are almost never used at capacity.

To evaluate the resilience of an ODMTS, this paper makes use of
\emph{scenarios}. Each scenario defines the vehicle capacities (e.g.,
to enforce social distancing policies), the frequencies, and the
demand for public transit.  Varying these inputs allows for evaluating
both pandemic and non-pandemic scenarios.
The scenarios focus on large shifts in demand due to the pandemic, rather than the daily variations which are very small in comparison.
The base scenario is the proposed 
ODMTS network designed for the estimated O-D matrix for the current system.
\changed{Improving the quality of service may attract new demand, which counteracts some of the decrease.
This effect is studied by \citet{BasciftciVanHentenryck2021-CapturingTravelMode}.
A similar analysis could be performed for the case study in this paper. However, this would require the development of new choice models to capture mode switching during a pandemic, which is a significant undertaking and is left for future research.}

\subsection{O-D Matrix Estimation}
\label{sec:method:odmatrix}

This paper uses techniques by \citet{BarryEtAl2002-OriginDestinationEstimation} to estimate the
O-D matrix from the transaction data and the APC data. Transactions
are first grouped into \emph{legs}, a leg being a part of a trip where
the rider is in a single bus or in the rail system.  After
establishing the individual legs, the legs are chained into trips, and
trip destinations are estimated when data is unavailable. This occurs because the transactions (obtained from the Breeze Cards) report \changed{tap-in} and \changed{tap-out} data for the trains, but only \changed{tap-in} data for the buses. However, since the ridership mostly consists of commuters, it is generally possible to infer the bus destination for a morning trip by using the corresponding bus boarding of the evening trip (and vice versa). The remaining \changed{destinations} can be estimated by using the APC data and ensuring that the distributions match under some reasonable assumptions. The result is a set of rider trips $T$, where each trip $t\in T$ is represented by an origin $o(t)$, a destination $d(t)$, and the number of riders $p(t)$.

\subsection{ODMTS Design}
\label{sec:method:design}
The network design uses a new optimization model that refines the model introduced by \citet{MaheoEtAl2019-BendersDecompositionDesign} to capture transfer constraints, multiple frequencies, and additional transportation modes.\footnote{A preliminary model was presented at the CPAIOR 2020 conference \citep{DalmeijerVanHentenryck2020-TransferExpandedGraphs}.}
The design problem chooses the fixed routes and their frequencies, minimizing a weighted combination of fixed and operation costs, and trip durations. Its solution implictly defines the legs taken by every rider.

\paragraph{Problem Description}

The design problem is defined in terms of a directed multigraph $G=(V,A)$ and a
set of trips $T$.  Vertices in $V$ represent geographical locations,
and arcs in $A$ are possible connections for the ODMTS design.  Arc
$a\in A$ captures the possibility of traveling from origin $i(a) \in
V$ to destination $j(a)\in V$ with mode $m(a) \in M$, where
$M=\{shuttle, bus, rail\}$.  Bus and rail arcs have a frequency $f(a)
\in F^{bus}$ and $f(a) \in F^{rail}$ respectively.  Bus and rail arcs
provide the connections between a set of \emph{transit hubs} $V_H
\subseteq V$, while shuttles provide the connections to and from the
hubs. Designing the network amounts to choosing which arcs $a \in A$
to open.
\changed{The shuttle arcs $a \in A$ capture that passengers are transported from $i(a)$ to $j(a)$ by shuttles, but the network design model does not specify how the shuttles perform this task.
Shuttle routes will be created dynamically as part of the simulation by the real-time shuttle dispatcher (Section~\ref{sec:method:rtdars}) to ensure a realistic evaluation of the system.}

This paper assumes that the rail network is fixed and that
shuttle connections can always be used, such that rail and shuttle
arcs are always opened. More precisely, the arc set $A$ is constructed
as follows.
Note that $G$ is a multigraph, and allows for parallel arcs corresponding to different modes and frequencies.
Rail arcs are added with a fixed frequency
between every pair of stations that is connected by a rail line, and it is
assumed that all rail stations are in $V_H$. Shuttle connections are
added for every trip: from the origin to the hubs, from the hubs to
the destination, and directly from the origin to the destination. Bus
arcs are added for all frequencies in $F^{bus}$ between all bus-only
hubs (hubs that are not rail stations) and between each bus-only hub
and the three nearest rail stations. The set of all bus arcs is
denoted by $A^{bus}$. What remains to decide is which bus connections
to provide and at which frequencies.

The objective of the design problem is to minimize a weighted
combination of \emph{fixed and variable costs} and \emph{trip
  duration}, with parameter $\alpha \in [0,1]$ used to balance these
objectives.  The cost is given the weight $1-\alpha$ and the total
trip duration, i.e., the total travel and waiting times, is given the
weight $\alpha$. The fixed cost for including a bus arc $a \in A$ is defined
as $\beta_a = (1-\alpha) \tau_a f_a c^{bus}$, i.e., the product of
travel time $\tau_a$, the number of buses over the time horizon $f_a$,
and hourly bus cost $c^{bus}$.  

The contribution to the objective of trip $t\in T$ traveling over arc
$a\in A$ is given by
\begin{equation}
    \gamma_a^t =
        \begin{cases}
            \mathrlap{p(t)\left((1-\alpha) d_a \bar{c}^{shuttle} + \alpha \tau_a\right)} \\
            & \textrm{ if } m(a) = shuttle\\
            p(t)\alpha\left(\tau_a + \frac{L}{2f(a)}\right) & \textrm{ if } m(a)\in \{bus,rail\}.
        \end{cases}
\end{equation}
The product $d_a \bar{c}^{shuttle}$ represents the variable cost of
using a shuttle along $a$, obtained by multiplying the distance by the
shuttle cost per mile.
This cost of using shuttle arc $a$ is an approximation for the actual cost of transporting a passenger from $i(a)$ to $j(a)$ when operating a real-time ride-sharing system (Section~\ref{sec:method:rtdars}).

The contribution to the trip
duration is $\tau_a$, since the shuttles are assumed to be available
within a short waiting time.  There are no variable costs associated
with using a bus or a train. The contribution to the trip duration is
the sum of the travel time and the expected waiting time.  This value
is calculated from the length of the time horizon $L$ and the number
of departures during that time.

Because transfers discourage people from using the ODMTS, a limit of
$K\ge 1$ legs is imposed on the length of each passenger path.  This
ensures that the total number of transfers is at most $K-1$, which
occurs when the passenger has to change vehicles between every pair of
arcs.
Vehicle capacity is not considered at the network design stage, motivated by the currently low utilization in Atlanta, but is still taken into account by the simulator.
Section~\ref{sec:baseline} shows that, for the case study, there is indeed no need to consider capacities, as passengers in the simulation almost always find an empty seat.

\paragraph{Optimization Model}

The decision variables of the model are specified as follows: Variable
$z_a \in \mathbb{B}$ is a binary variable with value one if and only if bus arc $a
\in A^{bus}$ is part of the design, while flow variable $y_a^t \in
\mathbb{B}$ indicates that the passengers of trip $t \in T$ travel
through arc $a \in A$. Formulation~\eqref{formulation:genodmts} presents the model.
For brevity, $\delta^+(i)$ denotes the set of
outgoing arcs for $i \in V$, and $\delta^+(i,m)$ the set of outgoing
arcs for $i\in V$ restricted to mode $m\in M$.  The sets $\delta^-(i)$
and $\delta^-(i,m)$ are defined similarly for incoming arcs.

\begin{mini!}
%
   	{}
%
   	{\sum_{a \in A^{bus}} \beta_a z_a + \sum_{t \in T} \sum_{{a} \in {A}} \gamma_a^t y_{{a}}^t, \label{eq:objectiveOveral}}
%
   	{\label{formulation:genodmts}}
%
   	{}
%
%
   	\addConstraint
   	{}
   	{\notag}
   	{\hspace{2.7 cm}\mathllap{\sum_{a \in \delta^+(i,bus)} f(a) z_a - \sum_{a \in \delta^-(i,bus)} f(a) z_a = 0}}
   	\addConstraint
   	{}
   	{\label{eq:masterFlowConservation}}
   	{\forall i \in V_H,}
   	\addConstraint
   	{\hspace{2.2 cm}\mathllap{\sum_{a \in A^{bus} \vert i(a) = i, j(a) = j} z_a}}
   	{\le 1 \quad \label{eq:masterOneFrequency}}
   	{\forall i, j \in V_H,}
    \addConstraint
    {}
    {\notag}
    {\hspace{2.5 cm}\qquad\mathllap{\sum_{{a} \in {\delta}^+({i})} y_{{a}}^t - \sum_{{a} \in {\delta}^-({i})} y_{{a}}^t = \begin{cases} 1 &\textrm{ if } {i} = o(t) \\ -1 &\textrm{ if } {i} = d(t) \\ 0 &\textrm{ else } \end{cases}}}
    \addConstraint
    {}
    { \label{eq:extsubFlowConservation}}
    {\forall t \in T, {i} \in {V},}
    \addConstraint
    {y_{{a}}^t}
    {\le z_a \label{eq:extsubArcCapacity}}
    {\forall t\in T, {a}  \in {A^{bus}},}
    \addConstraint
    {\sum_{a\in A} y_a^t}
    {\le K \label{eq:extsubTransferLimit}}
    {\forall t \in T,}
    \addConstraint
    {y_{{a}}^t}
    {\in \mathbb{B} \label{eq:extsubVariables}}
    {\forall t \in T, {a} \in {A},}
   	\addConstraint
   	{ z_a}
   	{\in \mathbb{B} \label{eq:masterVariables}}
   	{\forall a \in A^{bus}.}
\end{mini!}%

Objective~\eqref{eq:objectiveOveral} minimizes the cost for selecting
bus arcs, and the total passenger objective.
Constraints~\eqref{eq:masterFlowConservation} balance the number of
incoming and outgoing buses at every hub, while
Constraints~\eqref{eq:masterOneFrequency} ensure that at most one bus
frequency is selected for every connection.  Standard flow
conservation constraints for the passenger trips are given by
Constraints~\eqref{eq:extsubFlowConservation}.  Passengers only use
bus arcs that are part of the design, as stated by
Constraints~\eqref{eq:extsubArcCapacity}.  The transfer limit
constraint is enforced by Constraints~\eqref{eq:extsubTransferLimit},
and Equations~\eqref{eq:extsubVariables}-\eqref{eq:masterVariables}
are the integrality requirements.

\paragraph{Benders Decomposition}
\citet{MaheoEtAl2019-BendersDecompositionDesign} apply Benders decomposition to exploit the structure of the ODMTS design problem, and the same decomposition is used in this paper.
The \emph{master problem} determines the design by selecting bus arcs, and the \emph{subproblem} routes the riders through the network.
Let $\Phi(z)$ be the optimal value for the variable cost of the objective value 
given the values of the design variables $z$.
The master problem is given by Formulation~\eqref{formulation:genodmts:master}, and Formulation~\eqref{formulation:genodmts:sub} presents the subproblem.
The master problem concerns the design variables, and the passenger cost resulting from the design is implicitly defined through $\Phi(z)$.
The continuous variable $\theta$ is introduced for notational convenience, and equals $\Phi(z)$ in any optimal solution.
For a given design $z$, the subproblem is separable, and the problem of routing the riders through the network can be solved for each trip independently, which is a significant computational advantage.

\begin{mini!}
	%
	{}
	%
	{\sum_{a \in A^{bus}} \beta_a z_a + \theta, \label{eq:objective:master}}
	%
	{\label{formulation:genodmts:master}}
	%
	{}
	%
	%
	\addConstraint
	{\mathclap{\qquad\qquad\quad\eqref{eq:masterFlowConservation}, \eqref{eq:masterOneFrequency}, \eqref{eq:masterVariables},}}
	{\notag}
	{}
	\addConstraint
	{\theta}
	{\ge \Phi(z),}
	{\label{eq:thetalink}}
	\addConstraint
	{\theta}
	{\ge 0.}
	{}
\end{mini!}%

\begin{mini!}
	%
	{}
	%
	{\sum_{t \in T} \sum_{{a} \in {A}} \gamma_a^t y_{{a}}^t, \label{eq:objective:sub}}
	%
	{\label{formulation:genodmts:sub}}
	%
	{\Phi(z) = }
	%
	%
	\addConstraint
	{\eqref{eq:extsubFlowConservation}-\eqref{eq:extsubVariables}.}
	{\notag}
	{}
\end{mini!}%

Benders approaches construct an outer-approximation of Constraint~\eqref{eq:thetalink}, which is the epigraph of $\Phi(z)$, by iteratively solving the master problem for the current approximation, and then solving the subproblem to obtain a supporting hyperplane at the current value of $z$.
\citet{MaheoEtAl2019-BendersDecompositionDesign} show for their model that the subproblem is totally unimodular, which justifies the approach.

The new model~\eqref{formulation:genodmts}, however, is complicated by the transfer limit constraint~\eqref{eq:extsubTransferLimit}, which destroys total unimodularity and results in $\Phi(z)$ being non-convex.
Handling non-convex value functions has been addressed in the literature \citep{LaporteEtAl2002-IntegerLShaped,CodatoFischetti2006-CombinatorialBendersCuts}, but significantly increases the computational complexity. 
To remedy this limitation, the next section introduces an alternative formulation for the subproblem that enforces the transfer limit without destroying total unimodularity, so that a standard Benders decomposition method can be applied.

\paragraph{Transfer-Expanded Graphs}
To maintain total unimodularity in the subproblem, this section reformulates the subproblem on a \emph{transfer-expanded graph}.
Transfer-expanded graphs are layered graphs for which each layer corresponds to a number of transfers.
The transfer limit constraint is encoded directly by limiting the number of layers in the graph.
Transfer-expanded graphs share some similarities with time-expanded networks, where layers are used to represent different periods of time.
A literature review on time-expanded graphs is given by \citet{BolandEtAl2017-ContinuousTimeService}, and a more general discussion on layered graph approaches is presented by \citet{GouveiaEtAl2019-LayeredGraphApproaches}.

\begin{figure}
\centering
\clipbox{0cm 0cm 1cm 0cm}{
\begin{tikzpicture}[>=stealth', shorten >=1pt, auto, node distance=2cm, scale=0.5]

	\tikzstyle{origin}=[draw=black, circle, minimum size=0.5cm]
	\tikzstyle{dest}=[draw=black, circle, minimum size=0.5cm]
	\tikzstyle{hub}=[draw=black, fill=red, minimum size=0.5cm]
	\tikzstyle{shuttle}=[draw=gray, decoration={markings, mark=at position 0.6 with {\arrow{>}}}, postaction={decorate}]
	\tikzstyle{hubarc}=[->, draw=black, line width=0.4mm]
	
	\node[origin] (o) [] {};
	\node[hub] (h2k1) [below of=o] {};
	\node[hub] (h1k1) [left of=h2k1] {};
	\node[hub] (h3k1) [right of=h2k1] {};
	\node[hub] (h1k2) [below of=h1k1] {};
	\node[hub] (h2k2) [below of=h2k1] {};
	\node[hub] (h3k2) [below of=h3k1] {};
	\node[hub] (h1k3) [below of=h1k2] {};
	\node[hub] (h2k3) [below of=h2k2] {};
	\node[hub] (h3k3) [below of=h3k2] {};
	\node[dest] (d)	  [below of=h2k3] {};
	
	\path 	(o) edge[shuttle] (h1k1)
			(o) edge[shuttle] (h2k1)
			(o) edge[shuttle] (h3k1);
			
	\path	(o) edge[shuttle, bend left=80, looseness=1.4] (d);
	
	\path	(h1k1) edge[shuttle] (d)
			(h1k2) edge[shuttle] (d)
			(h1k3) edge[shuttle] (d)
			(h2k1) edge[shuttle, bend right=18] (d)
			(h2k2) edge[shuttle, bend right=20] (d)
			(h2k3) edge[shuttle] (d)
			(h3k1) edge[shuttle] (d)
			(h3k2) edge[shuttle] (d)
			(h3k3) edge[shuttle] (d);
			
	\path	(h1k1) edge[hubarc] (h2k2)
			(h2k1) edge[hubarc] (h3k2)
			(h1k2) edge[hubarc] (h2k3)
			(h2k2) edge[hubarc] (h3k3);
			
	\node[origin] (o) [label={above:$o(t)$}] {};
	\node[hub] (h2k1) [below of=o, label={above:$(2, 1)$}] {};
	\node[hub] (h1k1) [left of=h2k1, label={above:$(1, 1)$}] {};
	\node[hub] (h3k1) [right of=h2k1, label={above:$(3, 1)$}] {};
	\node[hub] (h1k2) [below of=h1k1, label={above:$(1,2)$}] {};
	\node[hub] (h2k2) [below of=h2k1, label={above:$(2,2)$}] {};
	\node[hub] (h3k2) [below of=h3k1, label={above:$(3,2)$}] {};
	\node[hub] (h1k3) [below of=h1k2, label={above:$(1,3)$}] {};
	\node[hub] (h2k3) [below of=h2k2, label={above:$(2,3)$}] {};
	\node[hub] (h3k3) [below of=h3k2, label={above:$(3,3)$}] {};
	\node[dest] (d)	  [below of=h2k3, label={below:$d(t)$}] {};
			
\end{tikzpicture}
}
\caption{Example \changed{of a} Transfer-Expanded Graph for $K=4$ and Two Hub Arcs in the Original Graph.}
\label{fig:teg}
\end{figure}
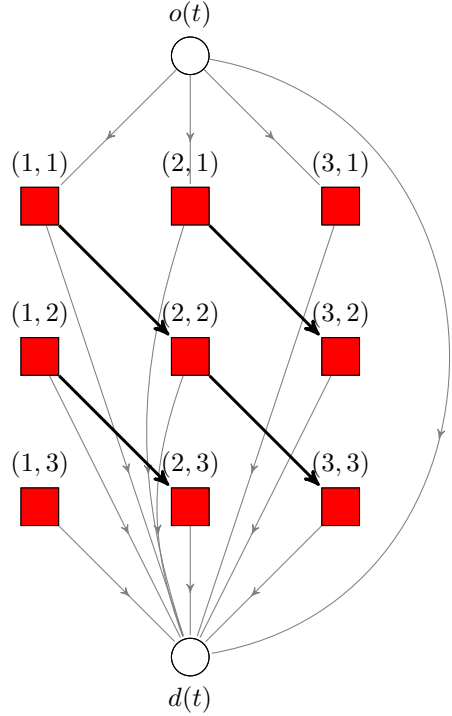

For clarity, the transfer-expanded graph is introduced for a given trip $t\in T$ (recall that the subproblem is separable).
Let $\bar{G}=(\bar{V}, \bar{A})$ denote the transfer-expanded graph.
The vertex set $\bar{V}$ consists of the origin $o(t)$, the destination $d(t)$, and $K-1$ copies of each hub, denoted by $(i, k)$ for $i \in V_H$ and $k \in \{1, \hdots, K-1\}$.
The arc set $\bar{A}$ is constructed as follows.
Figure~\ref{fig:teg} provides an illustrative example.
\begin{itemize}
	\item For every arc $a \in A$ between the hubs, add $K-2$ copies to $\bar{A}$, each going from one layer to the next.
	That is, add $\bar{a} = ((i(a), k), (j(a), k+1))$ for every $k \in \{1,\hdots,K-2\}$.
	\item Add $\bar{a} = (o(t), d(t))$ for the direct shuttle connection.
	\item Add shuttle arcs from the origin to the first layer of hubs, $\bar{a} = (o(t), (i, 1))$ for $i\in V_H$.
	\item Add shuttle arcs from every layer of hubs to the destination, $\bar{a} = ((i,k), d(t))$ for $i \in V_H$, $k\in \{1,\hdots,K-1\}$.
\end{itemize}

By construction, there is a bijection between paths in the transfer-expanded graph and paths in the original graph that satisfy the transfer limit constraint.
This allows the subproblem to be formulated as a standard shortest path problem on the transfer-expanded graph.
With a slight abuse of notation, let the cost $\gamma_{\bar{a}}^t$ of arc $\bar{a} \in \bar{A}$ be equal to $\gamma_a^t$ for the associated arc $a \in A$ in the original graph.
The binary variables $y_{\bar{a}}^t$ indicate the passenger flow through the transfer-expanded graph.
For notational convenience, let $\bar{A}(a)$ denote the copies of arc $a \in A$ in the transfer-expanded graph, and let $\bar{\delta}^+(\bar{v})$ and $\bar{\delta}^-(\bar{v})$ denote the outgoing and incoming arcs, respectively, for $\bar{v} \in \bar{V}$. Model~\eqref{formulation:genodmts:sub2} presents the new formulation for the subproblem.
Objective~\eqref{eq:objective:sub2} calculates the passenger cost in the same way as the original subproblem, and Constraints~\eqref{eq:extsubFlowConservation:sub2} are the equivalent of the flow conservation constraints~\eqref{eq:extsubFlowConservation}.
Constraints~\eqref{eq:extsubArcCapacity:sub2} ensure that copies of bus arcs in the transfer-expanded graph can only be used if the original arc is part of the design.
Finally, Equations~\eqref{eq:extsubVariables:sub2} define the passenger flow variables.
It is important to observe that the original transfer limit constraint~\eqref{eq:extsubTransferLimit} is encoded directly in the structure of the graph, and no additional constraint is needed.
The resulting subproblem is a minimum-cost flow problem, which has a totally unimodular constraint matrix, such that the integrality conditions~\eqref{eq:extsubVariables:sub2} may be relaxed.

\begin{mini!}
	%
	{}
	%
	{\sum_{\bar{a} \in \bar{A}} \gamma_{\bar{a}}^t y_{\bar{a}}^t, \label{eq:objective:sub2}}
	%
	{\label{formulation:genodmts:sub2}}
	%
	{}
	%
	%
	\addConstraint
	{}
	{\notag}
	{\hspace{2.7 cm}\qquad\qquad\qquad\mathllap{\sum_{\bar{a} \in {\bar{\delta}}^+(\bar{v})} y_{\bar{a}}^t - \sum_{\bar{a} \in {\bar{\delta}}^-(\bar{v})} y_{\bar{a}}^t = \begin{cases} 1 &\textrm{ if } {\bar{v}} = o(t) \\ -1 &\textrm{ if } {\bar{v}} = d(t) \\ 0 &\textrm{ else } \end{cases}}}
	\addConstraint
	{}
	{ \label{eq:extsubFlowConservation:sub2}}
	{ \forall {\bar{v}} \in {V},}
	\addConstraint
	{y_{\bar{a}}^t}
	{\le z_a \label{eq:extsubArcCapacity:sub2} \quad}
	{\forall {a}  \in {A^{bus}}, \bar{a} \in \bar{A}(a),}
	\addConstraint
	{y_{\bar{a}}^t}
	{\in \mathbb{B} \label{eq:extsubVariables:sub2}}
	{\forall \bar{a} \in \bar{A}.}
\end{mini!}%

The Benders decomposition method is applied as follows.
Let $\bar{\Phi}(z)$ be the optimal objective value of the subproblem defined on transfer-expanded graphs, which replaces $\Phi(z)$ in the master problem.
The constraint $\theta \ge \bar{\Phi}(z)$ is enforced by adding cutting planes.
For every $\hat{z} \in [0,1]^{\lvert A^{bus} \rvert}$, the standard \emph{Benders cut} \citep{Benders1962-PartitioningProceduresSolving} is given by
\begin{equation}
	\theta \ge \bar{\Phi}(\hat{z}) + \sum_{t \in T} \sum_{a \in A} \left(\sum_{\bar{a} \in \bar{A}(a)} \mu_{\bar{a}}^t(\hat{z})\right) (z_a - \hat{z}_a),
\end{equation}
where $\mu_{\bar{a}}^t(\hat{z})$ are optimal dual values for Constraints~\eqref{eq:extsubArcCapacity:sub2}, obtained by solving the subproblem for trip $t$ for design $\hat{z}$.
The network design problem is solved by repeatedly solving the current master problem, solving the independent subproblems to generate a Benders cut, and adding the cut to the master problem.
When no more cuts are violated, the problem is solved to optimality.
In this paper, the master problem is solved with CPLEX, and Benders cuts are separated in callbacks for both fractional and integer solutions $\hat{z}$.
The subproblems are solved with the LEMON graph library \citep{DezsoEtAl2011-LemonOpenSource}.

\subsection{Ridesharing and Fleet Sizing} \label{sec:method:fleetsizing}
Once the network is designed, the next step in the optimization
pipeline is to derive the appropriate size for the shuttle fleet.
This paper follows the same approach to fleet sizing as \citet{AuadVanHentenryck2021-RidesharingFleetSizing}, who study the impact of ridesharing and fleetsizing on ODMTS more generally.
The fleet sizing proceeds in two steps: it first identifies ridesharing
opportunities using a ridesharing optimization, before applying the
fleet-sizing optimization. The ridesharing optimization constructs
minimum-cost shuttle routes that cover the passenger shuttle requests,
and the fleet-sizing optimization finds the minimum number of shuttles
to cover these routes.  Both algorithms use the extended time horizon
$L^+$ to prevent that the fleet size is underestimated near the end of
the time horizon.

\paragraph{The Ridesharing Optimization}

Let $P$ be a set of shuttle requests, where each request $p \in P$
belongs to a single passenger and consists of an origin, a
destination, and a request time. The set $P$ is obtained from the
shortest paths through the designed network of all trips. It contains
a request for every shuttle arc used by riders.  If the request is at
the start of the trip, then the starting time of the trip is used as
the request time.  If this is not the case, the request time is
obtained by following the path. The shuttle requests are partitioned
based on the request type: for every hub $h \in V_H$, $P^-_h$ is
the set of requests with destination $h$ and $P^+_h$ the set of
requests with origin hub $h$.  The remaining requests, that
correspond to direct shuttle trips, are denoted by $P_{direct}$.

The objective of the ridesharing optimization is to construct shuttle
routes that cover all requests and minimize the total cost.
Riders are eligible to share the same shuttle if (i) the number of
passengers does not exceed the vehicle capacity $Q^{shuttle}$; (ii)
their departure times fall within a time window of width $\Delta \ge
0$; (iii) all the passenger either travel towards the same hub or
depart from the same hub (i.e., requests are in the same set $P^-_h$
or $P^+_h)$; and (iv) the duration of their shared route does not
exceed the duration of a direct trip by a factor more than tolerance
$\rho \ge 1$.

The route enumeration algorithm
by \citet{AuadVanHentenryck2021-RidesharingFleetSizing} is used to
construct the sets $R^-_h$ and $R^+_h$ of feasible routes for every
hub $h \in V_H$.  Routes in sets $R^-_h$ and $R^+_h$ are based on
requests from $P^-_h$ and $P^+_h$ respectively.  A route starts when
an empty shuttle picks up passengers, and ends when the last passenger
is dropped off.  Every shuttle route is associated with a starting
time and an ending time.  The starting time is the latest request time
among the passengers on the route, which ensures that the route can
only be started when all passengers are ready.  The ending time
follows from adding the route duration to the starting time.  The cost
of shuttle route $r$ is given by the parameter $\gamma^r$, which is
the same weighted sum of shuttle cost and passenger trip duration used
for designing the ODMTS.

The ridesharing optimization is solved as independent set partitioning
problems for every route set $R^-_h$ and $R^+_h$ associated with every
hub.  Without loss of generality, the formulation is presented for the
route set $R^-_h$ and the corresponding requests $P^-_h$.  Let
$x_r \in \mathbb{B}$ be a binary variable that takes the value of one
if route $r\in R^-_h$ is selected, and zero otherwise.  The binary
parameter $b^r_p$ is set to one if route $r\in R^-_h$ contains request
$p\in P^-_h$, and zero otherwise.

\begin{mini!}
	%
	{}
	%
	{\sum_{r \in R^-_h} \gamma^r x_r, \label{eq:rs:objective}}
	%
	{\label{formulation:rs}}
	%
	{}
	%
	%
	\addConstraint
	{\sum_{r\in R^-_h} b^r_p x_r}
	{= 1\quad \label{eq:rs:visitonce}}
	{\forall p\in P^-_h,}
	\addConstraint
	{x_r}
	{\in \mathbb{B} \label{eq:rs:vars}}
	{\forall r \in R^-_h.}
\end{mini!}%

Model~\eqref{formulation:rs} presents a formulation for
the ridesharing problem for route set $R^-_h$.
Objective~\eqref{eq:rs:objective} minimizes the total cost of the
selected routes.  Constraints~\eqref{eq:rs:visitonce} ensure that
every request is fulfilled by one of the routes.  Integral route
selections are enforced by Constraints~\eqref{eq:rs:vars}.
Problem~\eqref{formulation:rs} can be seen as a special case of the
problem solved
by \citet{AuadVanHentenryck2021-RidesharingFleetSizing}.

\paragraph{The Fleet-Sizing Optimization}

Given the optimal shuttle routes from Model~\eqref{formulation:rs},
and the single passenger routes related to $P_{direct}$, the
fleet-sizing algorithm finds the minimum number of shuttles necessary
to cover these routes. The fleet-sizing model uses a graph
$\mathcal{G} = (\mathcal{V}, \mathcal{A})$, where a vertex $r$
corresponds to a route, and an arc $(r, r')$ indicates that it is
feasible for a shuttle to serve route $r'$ after serving route $r$,
i.e., the ending time of route $r$, plus the time to relocate to the
start of route $r'$, is before the starting time of route $r'$.  The
graph $\mathcal{G}$ is extended with a source and a sink, denoted by
$src$ and $snk$ respectively.  The arc set $\mathcal{A}$ is extended
with arcs from the source to all other vertices and arcs from all
other vertices to the sink.

\begin{mini!}
	%
	{}
	%
	{\sum_{a \in \delta^+(src)} f_a, \label{eq:fs:objective}}
	%
	{\label{formulation:fs}}
	%
	{}
	%
	%
	\addConstraint
	{\sum_{a \in \delta^-(r)} f_a}
	{= \sum_{a \in \delta^+(r)}f_a \quad \label{eq:fs:flowconservation}}
	{\forall r \in \mathcal{V}\backslash\{src,snk\},}
	\addConstraint
	{\sum_{a \in \delta^-(r)} f_a}
	{= 1 \label{eq:fs:covering}}
	{\forall r \in \mathcal{V}\backslash\{src,snk\},}
	\addConstraint
	{f_a}
	{\ge 0 \label{eq:fs:vars}}
	{\forall a \in \mathcal{A}.}
\end{mini!}%

Formulation~\eqref{formulation:fs} presents the
fleet-sizing problem as a \changed{minimum flow problem with covering constraints}.
Let $f_a$ be the flow over arc $a \in \mathcal{A}$, and define
$\delta^-(r)$ and $\delta^+(r)$ to be the sets of arcs going into and
out of vertex $r\in \mathcal{V}$ respectively. By definition, a unit
flow from the source to the sink corresponds to a feasible schedule
for a single shuttle.  Therefore, Objective~\eqref{eq:fs:objective}
minimizes the size of the fleet by minimizing the total flow coming
out of the source.  Constraints~\eqref{eq:fs:flowconservation}
and \eqref{eq:fs:covering} ensure that the flows are balanced and that
every route is served.  Note that continuous variables may be used,
because the \changed{coefficient matrix is totally unimodular}.

\subsection{Real-Time Shuttle Dispatching}
\label{sec:method:rtdars}
To perform a realistic simulation of an ODMTS, it is necessary to
consider how the system is operated in real time and, in particular,
how to dispatch and route the shuttles.  To this end, the Real-Time
Dial-A-Ride System (RTDARS) by
\citet{RileyEtAl2019-ColumnGenerationReal} is embedded into the
simulator to dynamically dispatch shuttle requests.
\changed{The simulator provides the RTDARS with the current status of the shuttles in the fleet, and with shuttle requests based on the scenario demand.
Requests are only passed on as they arrive in real time, to prevent unrealistic anticipation.}

The RTDARS divides time into epochs of length $l$ and performs two
tasks for each period: it batches arriving shuttle requests during the
epoch, and it solves an optimization problem to route and dispatch
unserved requests from previous epochs.  A high-level overview of the
algorithm is presented here, and a more detailed description is
provided by \citet{RileyEtAl2019-ColumnGenerationReal}. The objective
of the static optimization problem at every epoch is to create minimum
cost shuttle routes that serve the requests.  These routes respect the
shuttle capacity $Q^{shuttle}$ and ensure that passengers do not
deviate too much from the shortest path.  Additionally, there is the
option to postpone requests to the next epoch by incurring a penalty.
The cost $c_r$ of route $r$ is the sum of the waiting times before
pickup.  The cost of not serving request $p$ is given by $g_p$, which
increases exponentially based on how long the rider has been waiting.
The initial penalty $\bar{g}$ doubles every ten epochs to ensure that
postponed requests are eventually served.

Let $P$ be the set of unserved shuttle requests from previous epochs,
let $R_s$ be the set of routes that can be assigned to shuttle $s \in
S$, and let $R = \cup_{s\in S} R_s$.  Feasible shuttle routes respect
previously committed actions, vehicle capacity, and a limit on the
time deviation from the shortest path.  The binary variable $\theta_r$
takes on the value one if and only if route $r \in R$ is part of the solution.
For every request $p\in P$, variable $w_p$ is equal to one if and only if the
request remains unserved in the current epoch.  The parameter $b^r_p$
indicates that request $p\in P$ is contained in route $r\in R$.
Formulation~\eqref{formulation:darp} presents the problem formulation.

\begin{mini!}
	%
	{}
	%
	{\sum_{r\in R} c_r \theta_r + \sum_{p \in P} g_p w_p, \label{eq:darp:objective}}
	%
	{\label{formulation:darp}}
	%
	{}
	%
	%
	\addConstraint
	{\sum_{r\in R} b^r_p \theta_r + w_p}
	{= 1 \quad \label{eq:darp:cover}}
	{\forall p \in P,}
	\addConstraint
	{\sum_{r\in R_s} \theta_r}
	{= 1 \label{eq:darp:assignment}}
	{\forall s \in S,}
	\addConstraint
	{\theta_r}
	{\in \mathbb{B} \label{eq:darp:routevars}}
	{\forall r \in R,}
	\addConstraint
	{w_p}
	{\in \mathbb{B} \label{eq:darp:skipvars}}
	{\forall p \in P.}
\end{mini!}%

Objective~\eqref{eq:darp:objective} minimizes the cost of the selected
routes and the penalties incurred for leaving requests unserved.
Constraints~\eqref{eq:darp:cover} specify that every request must be
served, or a penalty is incurred.
Constraints~\eqref{eq:darp:assignment} ensure that every shuttle is
assigned a single route.  Equations~\eqref{eq:darp:routevars} and
\eqref{eq:darp:skipvars} are the integrality
conditions. \citet{RileyEtAl2019-ColumnGenerationReal} solve the
static optimization problem with column generation.  This algorithm
iteratively generates promising shuttle routes (columns) with disjoint
sets of requests.  After reaching a stopping condition (based on
real-time requirements), Problem~\eqref{formulation:darp} is solved
for a restricted route set that only contain the routes that have been
generated.

\subsection{Simulation}
\label{sec:method:simulation}

The RTDARS is embedded in the simulator of
\citet{RileyEtAl2019-ColumnGenerationReal} to evaluate the performance
of the designed ODMTS.  The network design is provided by the design
problem, the number of shuttles is provided by the ridesharing
and fleet-sizing algorithms, and the demand is defined by the
scenario.  Passenger trips are revealed over time, and at every time
increment the simulator updates the position of every passenger and
every vehicle.  Results are obtained for time horizon $L$ by
simulating the extended horizon $L^+$ and truncating the results.

The original simulator only handled shuttle trips, and the following changes are made to simulate the ODMTS.
For every incoming trip request, the optimal path through the network
is determined as defined in Section~\ref{sec:method:design}.  When a
shuttle connection is encountered, a shuttle request is posted to the
RTDARS, while for bus and rail connections, passengers take the first
vehicle with an empty seat along their predetermined route.  The
vehicle capacities are given by $Q^{shuttle}$, $Q^{bus}$, and
$Q^{rail}$.  At the start of the time horizon, the shuttles are
distributed over the hubs proportionally to the average demand in each
area. The bus arcs in the design are grouped into lines, and each line
is assigned the minimum number of buses necessary to obtain the chosen
frequency.  Buses are spaced out evenly, but the rail and bus
schedules are not synchronized.  Solving a scheduling problem may
improve transfer times, and is an interesting direction for future
work. If the results show that the RTDARS is overwhelmed and shuttle
waiting times are long, then the fleet size is increased by 10\% and
the simulation is restarted.  This increase is to compensate for the
fact that the ridesharing and fleet-sizing optimizations assume
perfect information, while the RTDARS receives requests dynamically.

The simulator allows for high-fidelity analyses of the ODMTS system.
The analyses are based on simulated movements of riders and vehicles,
not on the design assumptions.  For example, the ODMTS design problem
uses expected waiting times, where simulation results allow for
calculating the actual waiting time for every rider.

\begin{table}[!t]
    \centering
    \scriptsize
    \begin{tabular}{lp{5.5 cm}}
        \toprule
        Parameter & Value \\
        \midrule
        $L$ & 4 hours (time horizon 6 a.m.-10 a.m.)\\
        $L^+$ & 6 hours (time horizon 6 a.m.-12 p.m.)\\
        $F^{bus}$ & $\{8, 12, 16\}$ (two, three, or four buses per hour)\\
        $F^{rail}$ & $\{24\}$ (six per hour)\\
        $\lvert V_H \rvert$ & 57 transit hubs (38 rail stations, 19 bus-only hubs)\\
        $\alpha$ & time is valued at \$7.25 per hour (US federal minimum wage)\\
        $\tau_a$ & shuttle and bus travel times from OpenStreetMap, rail from schedule\\
        $d_a$ & distances from OpenStreetMap\\
        $c^{bus}$ & \$72.15 per hour\\
        $c^{shuttle}$ & \$27.31 per hour\\
        $\bar{c}^{shuttle}$ & \$1 per mile\\
        $K$ & 4 arcs (at most three transfers)\\
        $\rho$ & tolerance factor 1.5\\
        $\Delta$ & 30 seconds\\
        $Q^{shuttle}$ & 4 passengers\\
        $Q^{bus}$ & 57 passengers\\
        $Q^{rail}$ & 576 passengers (192 for green line)\\
        $l$ & 30 seconds per epoch\\
        $\bar{g}$ & 420 seconds\\
        \bottomrule
    \end{tabular}
    \caption{The Parameter Values for the Case Study.}%
    \label{tab:parameter_values}%
\end{table}

\section{Baseline ODMTS for Atlanta}
\label{sec:baseline}

This section presents a baseline ODMTS for the city of Atlanta to
assess the resiliency of ODMTS during a pandemic. To provide the
proper context, the baseline ODMTS is compared to the existing transit
system in terms of service quality and cost.

\paragraph{Settings}

Table~\ref{tab:parameter_values} presents the parameter values used in
the case study, using the 2018 data introduced in
Section~\ref{sec:reviewmarta:data}. Recall that the ODMTS is designed
for the 6am--10am morning peak. The current rail system is fixed, with
a frequency of six per hour to approximate the MARTA schedule closely.
The 19 bus-only hubs are chosen by iteratively selecting a location
that is at least four miles from the other hubs, and has the most
passenger activity. Road travel times and distances were obtained from
OpenStreetMap using the GraphHopper library \citep{OpenStreetMap2020}.
Traffic congestion is not taken into account. Note however cities like
San Francisco have dedicated lanes for transit vehicles. Moreover, 
shuttle connections are typically local and short, and therefore less
affected by traffic.  For practical reasons, stops within a 1500 feet
radius are clustered, and 15,479 unique trips remain.

\paragraph{Network Design}

\begin{figure*}[!t]
    \centering
    \begin{minipage}{0.45\textwidth}
        \centering
        \includegraphics[trim=22cm 1cm 22cm 0cm,clip, height=0.4\textheight]{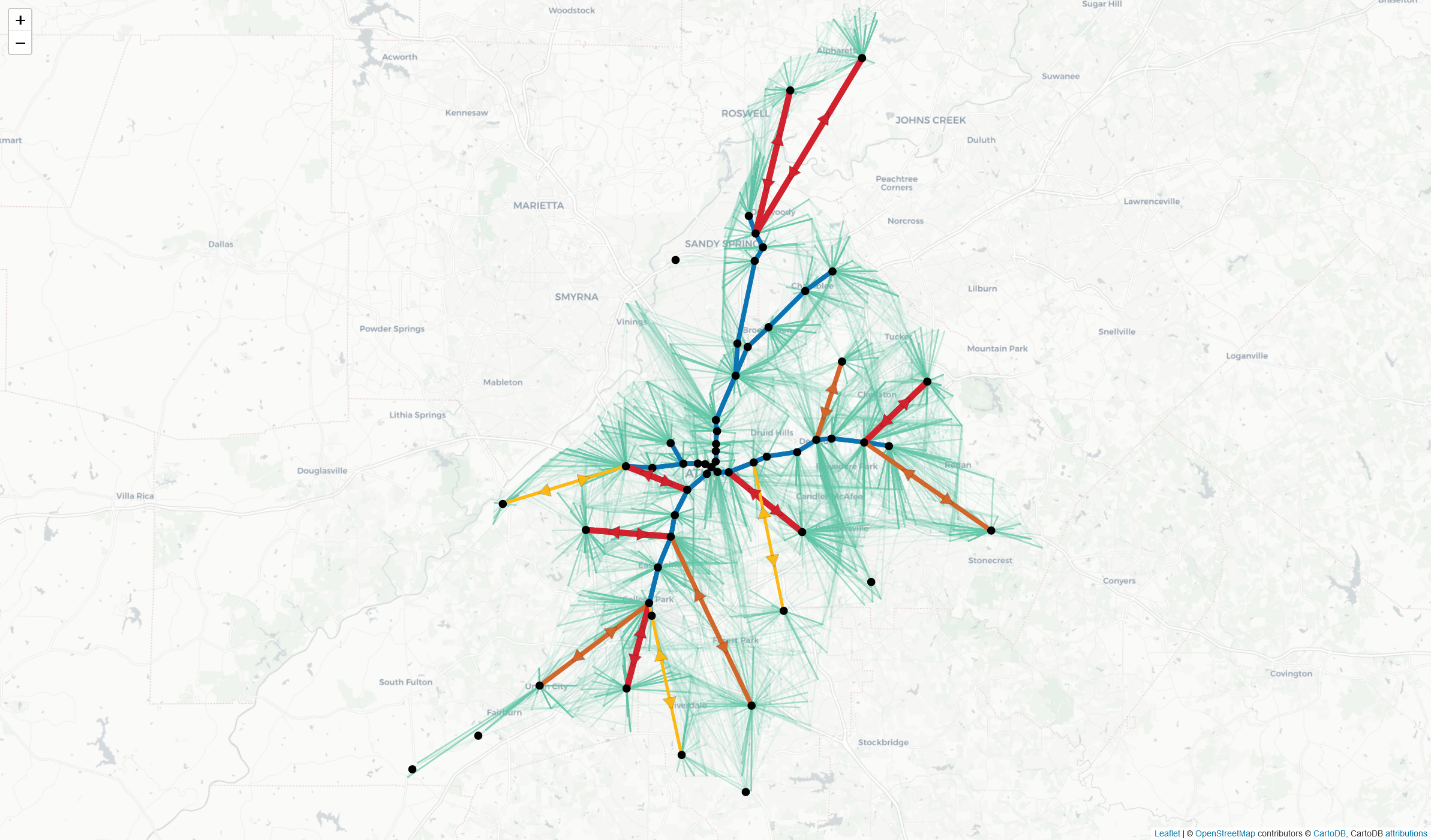}
        \caption{The Baseline ODMTS for Atlanta.}
        \label{fig:design}
    \end{minipage}%
    \begin{minipage}{0.1\textwidth}
        ~
    \end{minipage}%
    \begin{minipage}{0.45\textwidth}
        \centering
        \includegraphics[height=0.4\textheight]{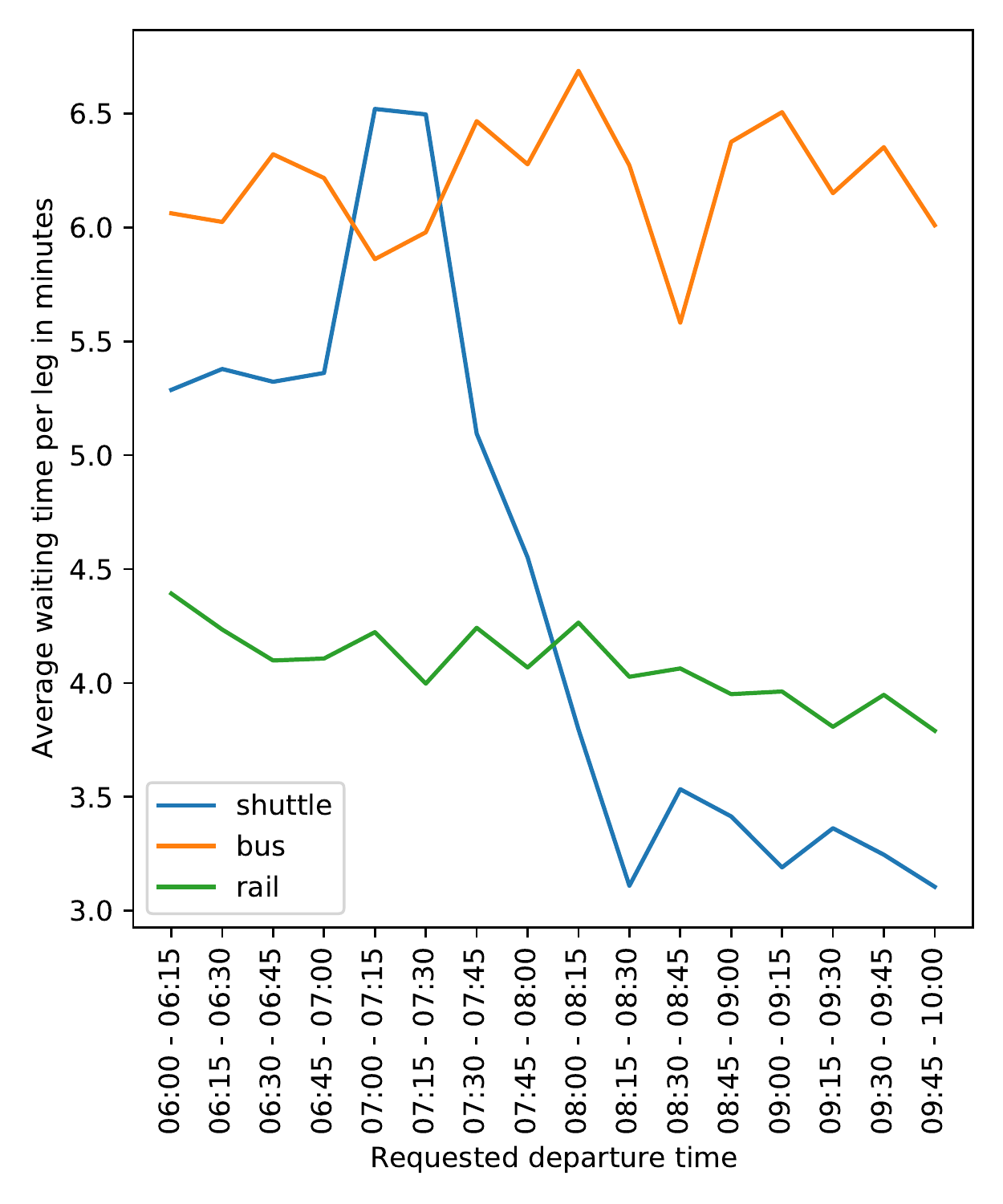}
        \caption{Average Waiting Time per Leg Type over Time.}
        \label{fig:baseline:time:waitingleg}
    \end{minipage}
    \\
    \vspace{1 cm}
    \begin{minipage}{0.45\textwidth}
        \centering
        \includegraphics[height=0.4\textheight]{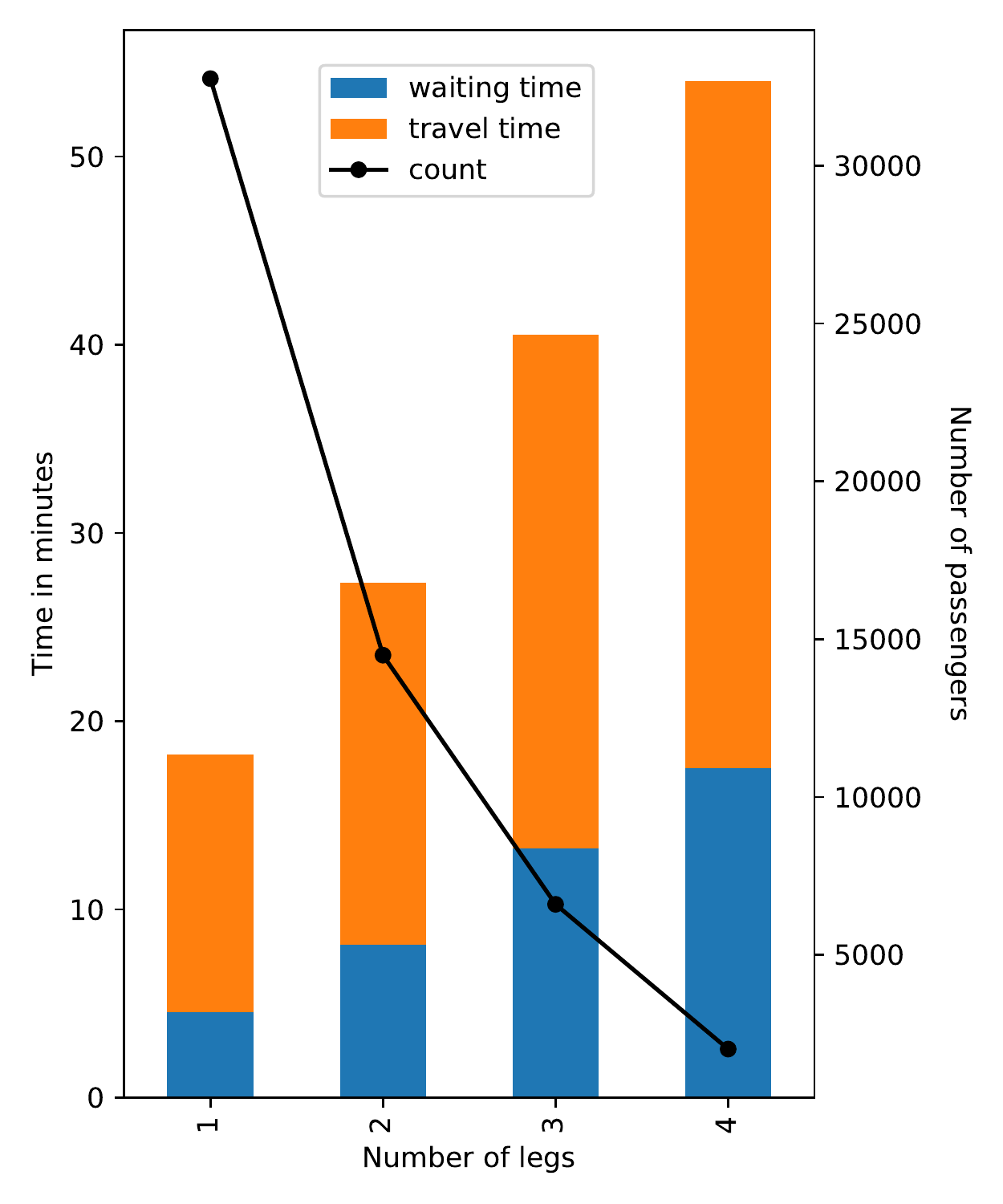}
        \caption{Average Passenger Wait and Travel Times (left axis) and Total Trip Counts (right axis) for a Given Number of Legs per Trip.}
        \label{fig:baseline:legs:combined}
    \end{minipage}%
    \begin{minipage}{0.1\textwidth}
        ~
    \end{minipage}%
    \begin{minipage}{0.45\textwidth}
        \centering
        \includegraphics[height=0.4\textheight]{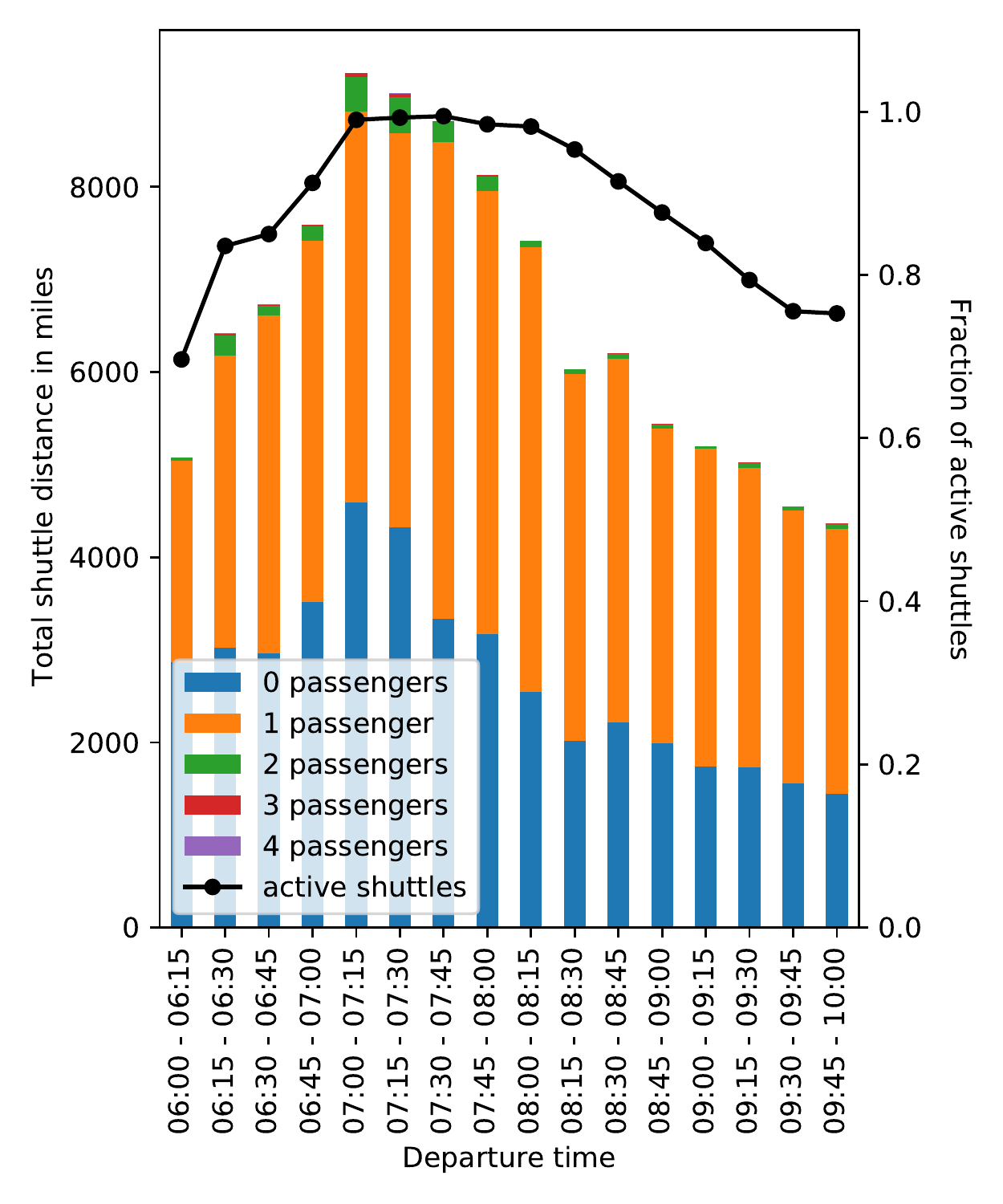}
        \caption{Total Miles Driven at Given Occupancy (left axis) and Fraction of Shuttles Active over Time (right axis).}
        \label{fig:baseline:capacity:combined}
    \end{minipage}
\end{figure*}

The resulting ODMTS design is shown in Figure~\ref{fig:design}.
Shuttle connections are depicted as thin lines (green), rail
connections as thick lines (blue), and bus connections as arrows
ranging from light to dark (yellow to red) for low to high
frequencies.  The rail is obviously the core of the network. Buses are
mostly used to expand the rail, providing connections to high-density
areas that it cannot reached. This holds for all four cardinal directions,
and especially for the North. Other bus lines cover the Northeast, the Southeast, and the
Southwest.  The Northwest is served by shuttles, because most
riders are close to the rail.  The ridesharing and fleet-sizing
optimizations estimate that 917 shuttles are necessary to operate the
ODMTS.  This number was scaled up by 20\% to 1100 shuttles to obtain
an acceptable service level. The bus service is provided by 24 buses.

\paragraph{Passenger Convenience}

\begin{table}[!t]
  \centering
    \begin{tabular}{p{1.5 cm}lrrr}
    \toprule
          &       & \multicolumn{3}{c}{Wait time (min)} \\
    \cmidrule{3-5}
          &       & \multicolumn{1}{r}{0-5} & \multicolumn{1}{r}{5-10} & \multicolumn{1}{r}{$>$10} \\
    \midrule
    \multirow{3}[0]{*}{Baseline} & Shuttle & 57\%  & 37\%  & 6\% \\
          & Bus   & 45\%  & 38\%  & 18\% \\
          & Rail  & 63\%  & 32\%  & 5\% \\
    \bottomrule
    \end{tabular}%
    \caption{Distribution of Total Waiting Times per Mode for Trips Using that Mode (Baseline).}
  \label{tab:baseline:passenger:waiting}%
\end{table}%

In real-time simulations, the average trip duration was found to be 25
minutes, with 18 minutes of travel and 7 minutes of waiting time.
Rail-only riders (56\%) spend 23 minutes in transit, which is similar to the estimate for the current system.  For the remaining trips (44\%), the trip
duration is only 27 minutes for the ODMTS compared to 46 minutes for
the current system.
Table~\ref{tab:baseline:passenger:waiting} shows the
distribution of the waiting times over the different modes. Only 6\%
of the trips including a shuttle leg requires waiting more than 10
minutes for a shuttle. The bus and rail distributions provide a
similar picture: waiting under five minutes is most common, and
waiting longer than 10 minutes is least likely. Shuttle waiting times vary
over time, while bus and rail waiting times are fairly constant.
Figure~\ref{fig:baseline:time:waitingleg} shows that the average
waiting time is about 6 minutes for a bus and 4 minutes for a train
leg. The shuttle waiting times have a peak of about 6.5 minutes
between 7am and 8am, which is the period that sees the most
passengers.  Around 6am, waiting times are also relatively long while
the shuttles spread out over the region.  After 8am, the shuttle
waiting time rapidly decreases to only 3 minutes.
Figure~\ref{fig:baseline:legs:combined} shows how waiting times are
distributed over the trips.  Short trips have short waiting times, and
waiting times increase almost proportionally with the length of the
trip, which indicates a fair distribution.  Furthermore, the counts
show that trips with many transfers are rare: the vast majority of
riders are served either by a direct trip or with a single transfer.

\paragraph{Capacity Utilization}

\begin{figure}[!t]
    \includegraphics[width=\linewidth]{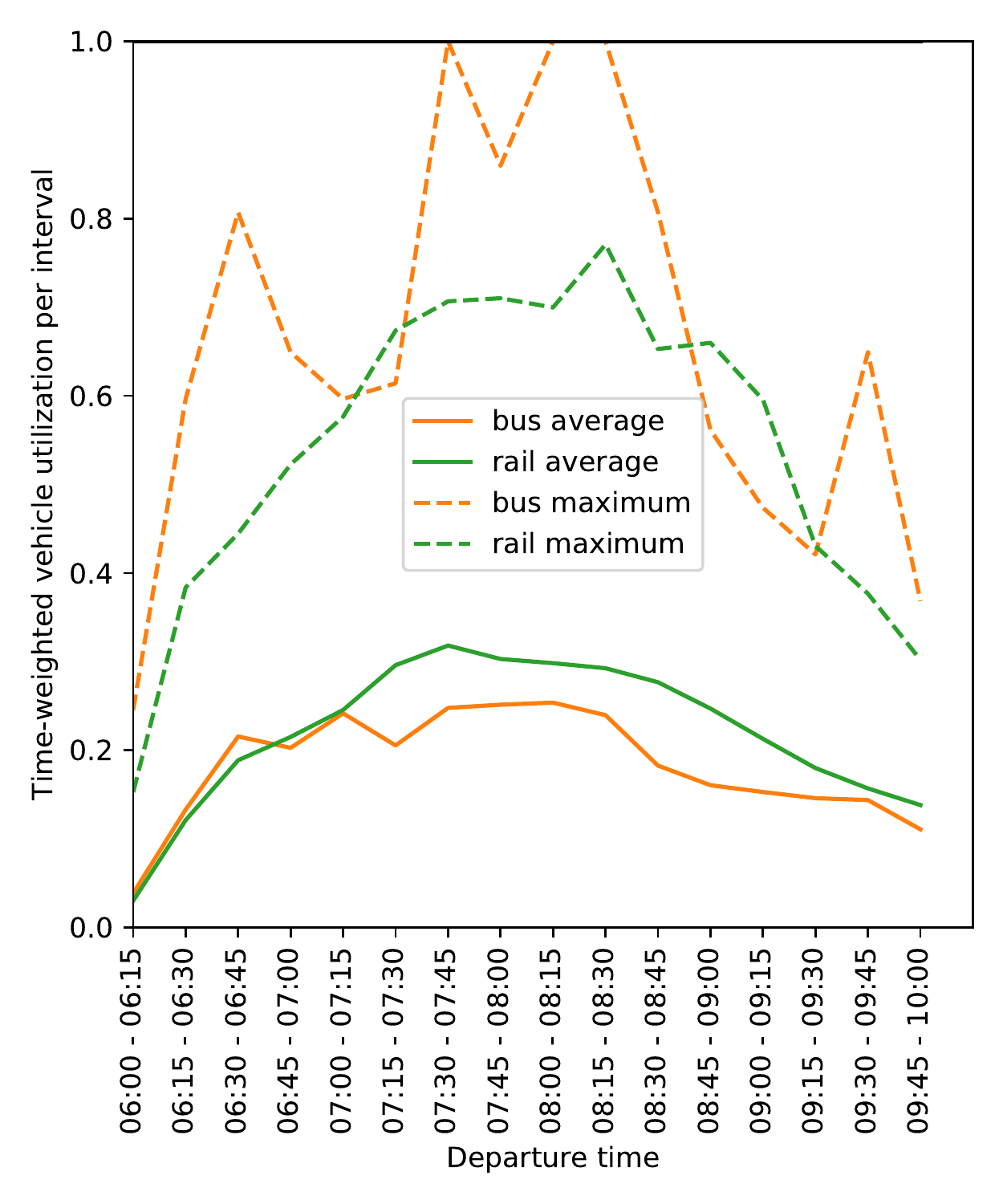}
    \captionof{figure}{Utilization of the Buses and Trains over Time (Baseline).}
    \label{fig:baseline:time:occupation}
\end{figure}

It is interesting to report capacity utilization results, since social
distancing measures will reduce capacity. The average utilization of
buses and rail stays under 32\% over the time horizon.
 The line chart in
Figure~\ref{fig:baseline:capacity:combined} shows the fraction of
active shuttles over the time horizon.  Almost all shuttles are active
during the busiest period from 7am to 8am, but the waiting times show
that the system is not overwhelmed. The same figure presents the
total mileage for a given shuttle occupancy.  It shows that riders
can be offered individual shuttle rides most of the time, and that the
amount of sharing increases when more shuttles are busy.  The low
occupancy suggests that reducing shuttle capacity in a pandemic is
feasible.  In total, the shuttles drive 105,127 miles between 6am and
10am, which amounts to 96 miles per shuttle on average. Trains are never at capacity, and only 0.7\% of bus passengers incur a delay due to buses being full. Figure~\ref{fig:baseline:time:occupation}
provides additional details.

\begin{figure*}[!t]
	\centering
	\begin{minipage}{0.5\linewidth}
		\centering
		\includegraphics[frame=1pt, width=0.95\linewidth, trim=45.5cm 4cm 41.5cm 0cm, clip]{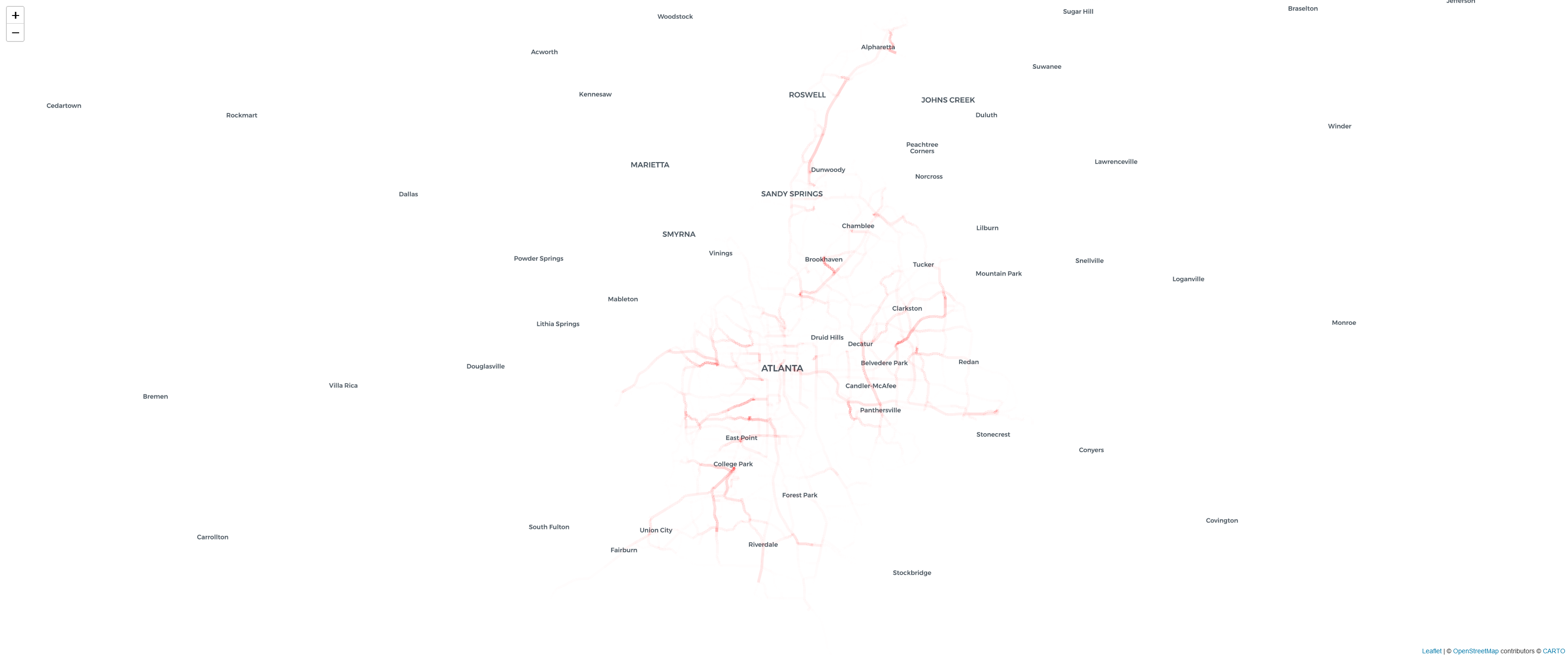}
		\caption{Road Usage with the ODMTS.}
		\label{fig:road_odmts}
	\end{minipage}%
	\hfill
	\begin{minipage}{0.5\linewidth}
		\centering
		\includegraphics[frame=1pt, width=0.95\linewidth, trim=45.5cm 4cm 41.5cm 0cm, clip]{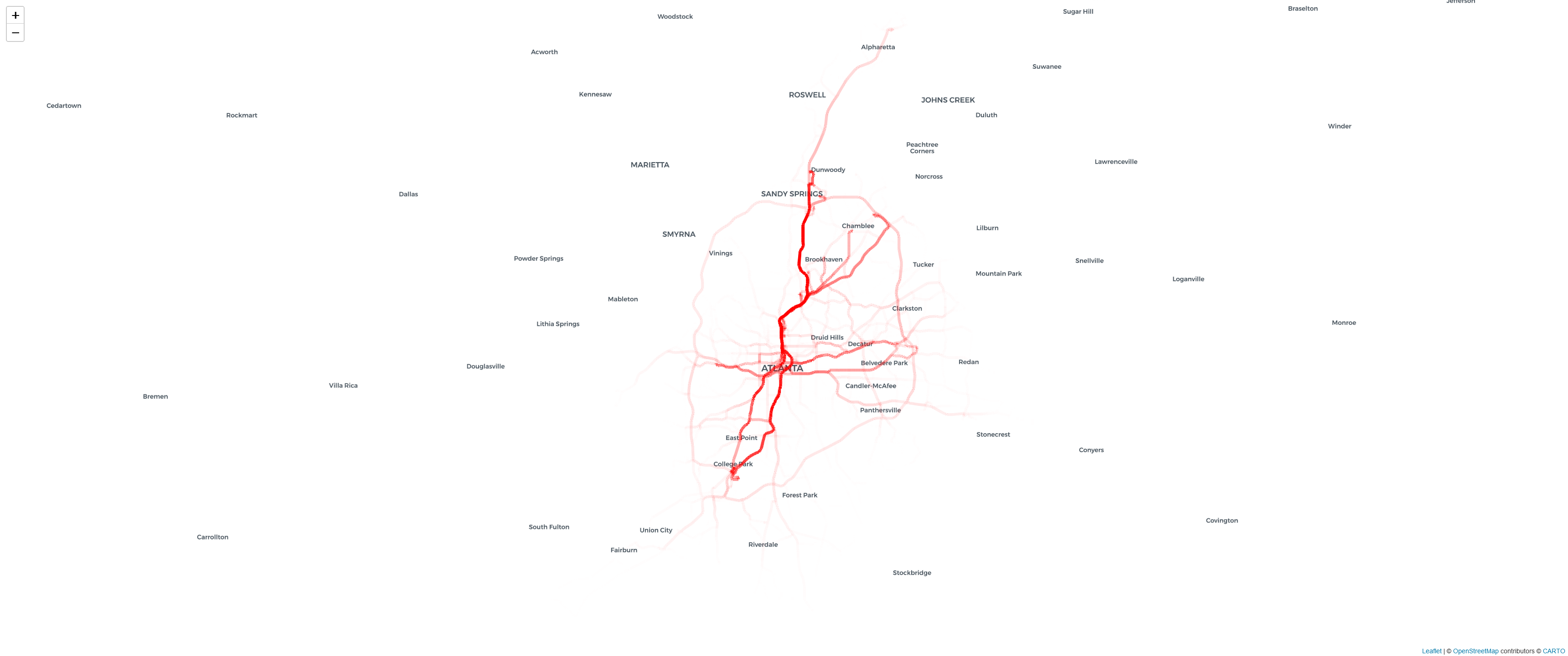}
		\caption{Road Usage Without the ODMTS.}
		\label{fig:road_car}
	\end{minipage}\\
	\vspace{0.5cm}
	\begin{minipage}{0.6\linewidth}
		\centering
		\includegraphics[frame=1pt, width=0.95\linewidth, trim=25cm 10cm 27.5cm 10cm,clip]{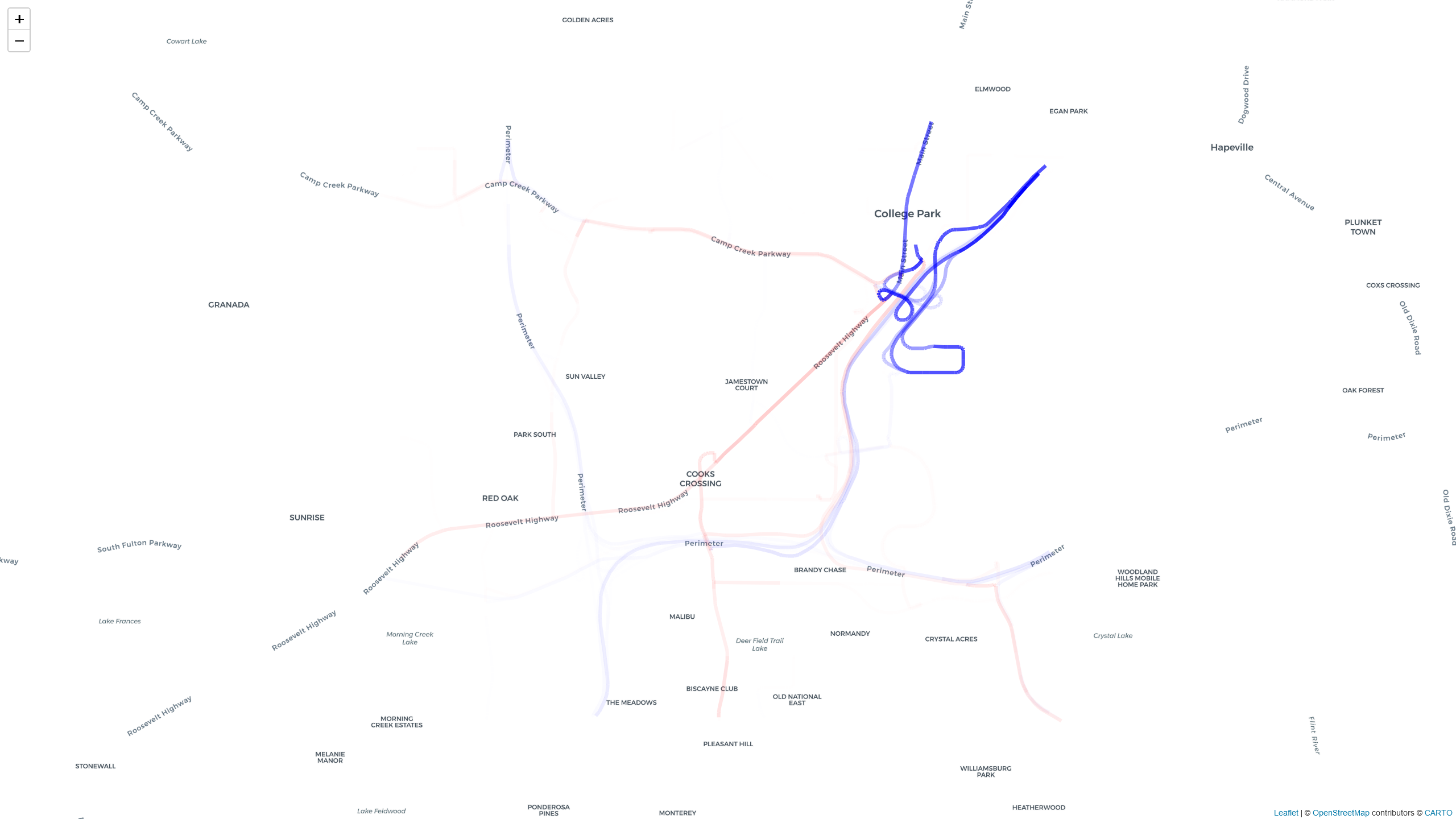}
		\caption{Change in Road Usage around College Park due to the ODMTS (blue for improvement).}
		\label{fig:road_diff}
	\end{minipage}
\end{figure*}

\paragraph{Road Usage and Congestion}

Figure~\ref{fig:road_odmts} visualizes the impact of the baseline ODMTS on the road network, where the opacity of each road increases with the number of times it is used by buses and shuttles. The figure shows that the ODMTS performs as intended, with mainly local shuttle trips. The shuttle trips impose little traffic pressure on the main roads, which are served by a relatively small number of buses.
Hence, by providing a convenient and accessible alternative to driving a car, ODMTS can also help to reduce congestion. For comparison, Figure~\ref{fig:road_car} visualizes the road usage \changed{when the same travelers would drive their personal vehicles.
The figure suggests that if more people would switch from car to ODMTS, this could significantly reduce traffic, and especially relieve congestion on the main roads, which are most affected.}
Some roads (for instance) close to the College Park train station may seem to have additional traffic. Figure \ref{fig:road_diff} \changed{provides a close-up of} this region. In fact, the ODMTS generates at most an average of 3.5 cars per minute more (red roads). The ODMTS also reduces the number of cars by up to 12.1 cars per minute \changed{compared to personal vehicles} (blue roads).

\paragraph{Comparison With the Current System}

To provide proper context, the baseline ODMTS is compared to the current system in terms of estimated cost and service quality.
While this paper focuses on the resiliency of ODMTS, a more detailed comparison between the two systems is an interesting direction for future research.
The estimated cost for buses in the morning peak is \$134k per day for the
current system. For the baseline ODMTS, it is only \$127k for
buses and shuttles together. These costs are calculated based on the
hourly vehicle costs.  Lower costs stem from the relatively low cost
of operating shuttles compared to buses: The ODMTS only uses buses for
busy corridors, where high-capacity vehicles are most efficient.  For
rail-only passengers, there is no significant difference between the
systems obviously, but it is not known how these riders travel to the rail
station. For current bus users, the ODMTS seems strongly preferred, as
the trip duration is 19 minutes shorter on average. Other advantages
of the ODMTS include better access to the public transit system, and
less time for riders to get to a stop, as the on-demand shuttles
can pick them up close to their origins. 

\paragraph{Importance of Buses}

As mentioned earlier, the core of the ODMTS for Atlanta is the rail
network, and bus lines provide access to the rail for high-density
areas.  Shuttles, however, can perform the same task. This gives rise
to the question of whether buses are important to the baseline ODMTS.
The baseline ODMTS has 24 buses and 1100 shuttles. If buses are
removed from the system, the ridesharing and fleet sizing algorithms
indicate that 151 additional shuttles are necessary, for a total fleet
of 1251 shuttles. The cost for the morning peak would then increase
by 7.5\% from \$127k to \$137k, which exceeds the cost of the
current system of \$134k. This is consistent with the case study of
ODMTS in Ann Arbor, where buses were shown critical
\citep{BasciftciVanHentenryck2020-BilevelOptimizationDemand}.

\section{Resiliency of ODMTS}
\label{sec:casestudy}

The resiliency of the ODMTS during a pandemic is evaluated with
scenarios that correspond to different periods of the COVID-19
pandemic. These scenarios differ in the level of ridership, and the
amount of social distancing measures. The resiliency is assessed by
investigating how the performance and cost evolve when
the baseline design is kept in place and only the number of shuttles
is scaled down. This is complemented by an evaluation that redesigns
the network to take into account the change in demand and additional cleaning costs.
Keeping the baseline ODMTS simplifies operations but a
redesign may yield additional cost savings. For MARTA, only 33\% of the operating expenses come from fare
revenues, but this number may be significantly higher for other
transit agencies. This section also studies how service levels change
when this percentage increases. 

\subsection{Scenarios}
\label{sec:casestudy:scenarios}

\begin{figure}[!t]
\centering
\includegraphics[trim=22cm 1cm 22cm 0cm,clip,width=0.9\linewidth]{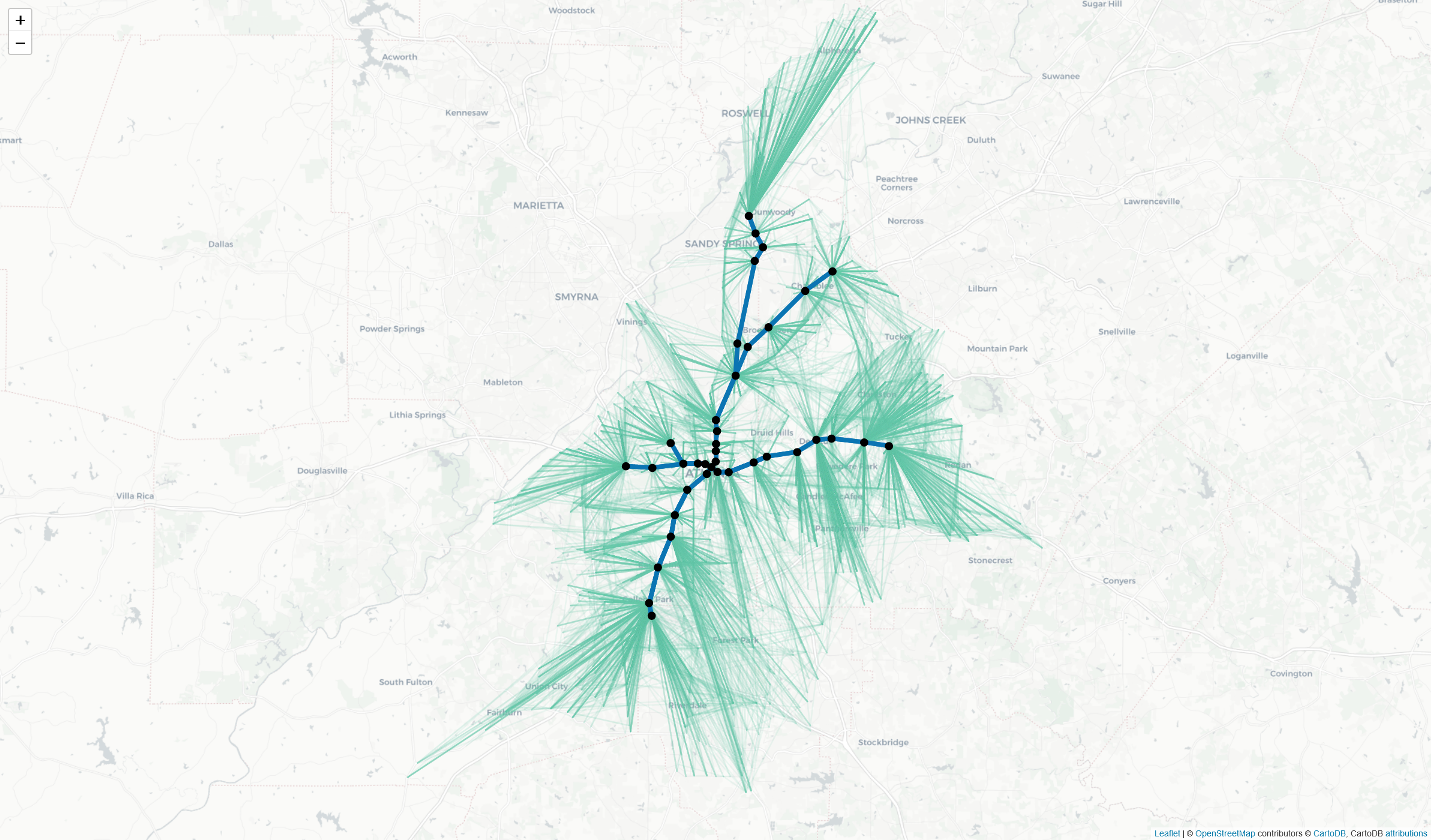}
    \captionof{figure}{The ODMTS for Atlanta Without Buses.}
    \label{fig:design_nobus}
\end{figure}

Table~\ref{tab:scenarios} summarizes the scenarios used in the case
study.  The first scenario is the baseline, and the other three are
pandemic scenarios with depressed demand and reduced vehicle capacities to
allow for social distancing.  Fleet sizes are provided by the
ridesharing and fleet-sizing algorithms.  The ``Fixed design'' column
assumes that the baseline network design is kept in place for all
scenarios, except for the strict late pandemic scenario, which removes
the bus subsystem (see Figure~\ref{fig:design_nobus}).  For the ``Redesign'' column, the network
is redesigned based on the scenario, and the fleet size is
determined for the redesigned network (Section~\ref{sec:case_study:redesign}).

\begin{table*}[!t]
  \footnotesize
  \centering
    \begin{threeparttable}
    \begin{tabular}{lllccccccc}
    \toprule
          &       &    &     \multicolumn{3}{c}{Capacity} && \multicolumn{2}{c}{Shuttle fleet}\\
    \cmidrule{4-6} \cmidrule{8-9}
    \multicolumn{1}{l}{Scenario name} & Ridership level & Passengers & Shuttle & Bus & Rail && Fixed design & Redesign & Budget\\
    \midrule
    Baseline & 2018 & 55871 (100\%) & 4 & 100\% & 100\% && 1100 &  & \$127k\\
    Early Pandemic & March 2020 & 25042 (45\%) & 1 & 50\%  & 50\%  && 846 & 837 & \$106k\\
    Late Pandemic & April 2020 & 13479 (24\%) & 1 & 50\%  & 50\%  && 510 & 514 & \$98k\\
    Strict Late Pandemic & April 2020 & 13479 (24\%) & 1 & 0\%  & 25\%  && 602\tnote{*} & 599 & \$98k\\
    \bottomrule
    \end{tabular}%
    \begin{tablenotes}
      \item[*] bus subsystem removed from the design.
    \end{tablenotes}
    \end{threeparttable}
    \caption{Overview of Scenarios Settings and Shuttle Fleets.}%
    \label{tab:scenarios}%
\end{table*}%

The passenger trips in the pandemic scenarios are sampled from the
2018 data until the ridership matches the pandemic level.  Sampling is
used because the trips observed during the pandemic are shaped by the
essential service plan, and do not necessarily represent the actual
demand.  The early pandemic scenario is based on the ridership level
around March 18th, 2020, with 70k bus and 60k rail transactions per
day.  The late pandemic scenarios use the April 2020 level with 40k
bus and 35k rail transactions. Due to the lack of bus transaction data
for April, the proportion of bus to rail trips is kept the same for
April and scaled down based on number of rail trips.  The bus and rail
capacities are indicated as percentages of the total baseline
capacity.  A 50\% capacity allows for only using every other seat for
social distancing.  As social distancing is difficult in a shuttle,
only one passenger is allowed during the pandemic.  The strict late
pandemic scenario imposes even stricter measures: Buses are considered
not safe, and the system operates with only single passenger shuttles
and trains at 25\% capacity.

One interesting observation is that for all pandemic scenarios, the
number of shuttles is smaller than that for the baseline scenario.  It
is not obvious that this would be the case, and it demonstrates that
ridership decreases sufficiently to allow for both a reduction in
capacity \emph{and} a reduction in the number of shuttles.  It follows
that the baseline fleet is sufficiently large to operate the system
during a pandemic, and no additional vehicles are necessary.

Every scenario is assigned a budget for buses and shuttles that takes
into account that reduced ridership affects fare revenue.  For the
current system, 33\% of the operating expenses are paid from fares and
direct revenue, while capital expenditures are completely funded from
other sources.  The baseline budget is taken to be the normal cost of
the system, 92\% of which is operating expenses (labor and
maintenance) and 8\% is capital expenditures (depreciation).  For the
other scenarios, the budget is reduced proportionally to the change in
available funds.  For the early pandemic scenario, for example, the
number of passengers is reduced by 55\%.  This results in a reduction
of operating funds by 18\%, and therefore a reduction of the budget by
17\% to \$106k.

\subsection{Resiliency under Fixed Design}
\label{sec:case_study:fixed}

This section evaluates the ODMTS resiliency when the baseline network
design stays in place during the
pandemic. Table~\ref{tab:pandemic_overview} presents the aggregate
results.  The average trip durations suggest that the ODMTS is very
resilient in terms of rider convenience. For riders not using rail,
the average trip duration remains under 27 minutes on average. These
results are also valid for individual trips.
Table~\ref{tab:passenger:waiting} shows that, for the early and late
pandemic scenarios, the wait times are comparable to the baseline in
distribution.  In the strict late pandemic scenario, waiting times for
the shuttles go up, as a result of removing the buses.

\begin{figure}[!t]
    \centering
    \includegraphics[width=0.8\linewidth]{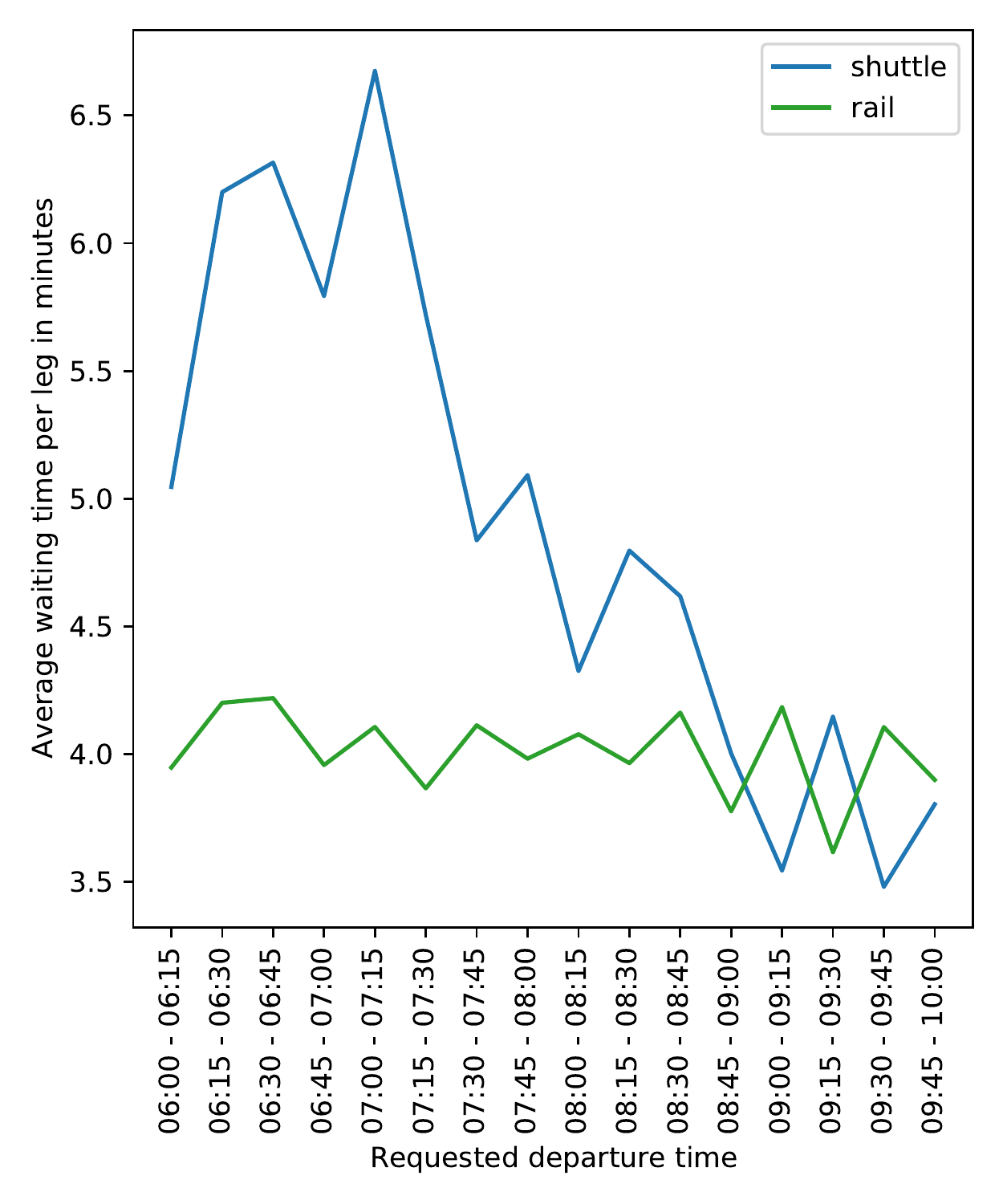}
    \captionof{figure}{Average Waiting Time over Time for a Leg of a Given Mode (Strict Late Pandemic).}
    \label{fig:strictlate:time:waitingleg}
\end{figure}

Figure~\ref{fig:strictlate:time:waitingleg} displays the waiting times
over time for the strict late pandemic scenario. The pattern is similar to the baseline,
except that the shuttle waiting times are slightly larger. But, even during the peak, the average
waiting time for a shuttle leg is at an acceptable 6.5 minutes.  This
shows that a high quality of service can be provided, even if buses
are no longer used, at the cost of increasing the number of
shuttles from 510 to 602.

\begin{table*}[!t]
  \centering
    \begin{tabular}{lrrrrrrrr}
    \toprule
          & \multicolumn{2}{c}{Average trip} && \multicolumn{2}{c}{Max. use of} &       &  \\
          & \multicolumn{2}{c}{duration (min)} && \multicolumn{2}{c}{total capacity} &       & \multicolumn{2}{c}{Cost (\$)}\\
          \cmidrule{2-3} \cmidrule{5-6} \cmidrule{8-9}
          & Rail-only & Other && Bus & Rail & Budget (\$) & Ex. Clean. & Incl. Clean. \\
          \midrule
    Baseline & 23    & 27    && 25\%  & 32\%  & 127k  & 127k  & 127k\\
    Early pandemic & 22    & 26    && 17\%  & 13\%  & 106k  & 99k & 105k\\
    Late pandemic & 22    & 26    && 10\%   & 7\%   & 98k & 63k  &  66k \\
    Strict late pandemic & 22    & 25    && 0\% & 7\%  & 98k & 66k & 70k \\
    \bottomrule
    \end{tabular}%
  \caption{Statistics per Scenario Under the Fixed Design.}
  \label{tab:pandemic_overview}%
\end{table*}%

\begin{table}[!t]
  \centering
    \begin{tabular}{p{1.5 cm}lrrr}
    \toprule
          &       & \multicolumn{3}{c}{Wait time (min)} \\
    \cmidrule{3-5}
          &       & \multicolumn{1}{r}{0-5} & \multicolumn{1}{r}{5-10} & \multicolumn{1}{r}{$>$10} \\
    \midrule
    \multirow{3}[0]{*}{Baseline} & Shuttle & 57\%  & 37\%  & 6\% \\
		& Bus   & 45\%  & 38\%  & 18\% \\
		& Rail  & 63\%  & 32\%  & 5\% \\
    \midrule
    \multirow{3}[0]{*}{\parbox{2 cm}{Early\\pandemic}}  & Shuttle & 62\%  & 35\%  & 4\% \\
          & Bus   & 44\%  & 36\%  & 20\% \\
          & Rail  & 64\%  & 31\%  & 5\% \\
    \midrule
    \multirow{3}[0]{*}{\parbox{2 cm}{Late\\pandemic}} & Shuttle & 62\%  & 35\%  & 4\% \\
          & Bus   & 44\%  & 39\%  & 17\% \\
          & Rail  & 64\%  & 31\%  & 5\% \\
    \midrule
    \multirow{3}[0]{*}{\parbox{2 cm}{Strict late\\pandemic}} & Shuttle & 53\%  & 37\%  & 10\% \\
          & Bus   & n.a.  & n.a.  & n.a. \\
          & Rail  & 61\%  & 33\%  & 6\% \\
    \bottomrule
    \end{tabular}%
    \caption{Distribution of Total Waiting Times per Mode for Trips using that Mode.}
  \label{tab:passenger:waiting}%
\end{table}%

\begin{figure}[!t]
    \centering
    \includegraphics[width=0.8\linewidth]{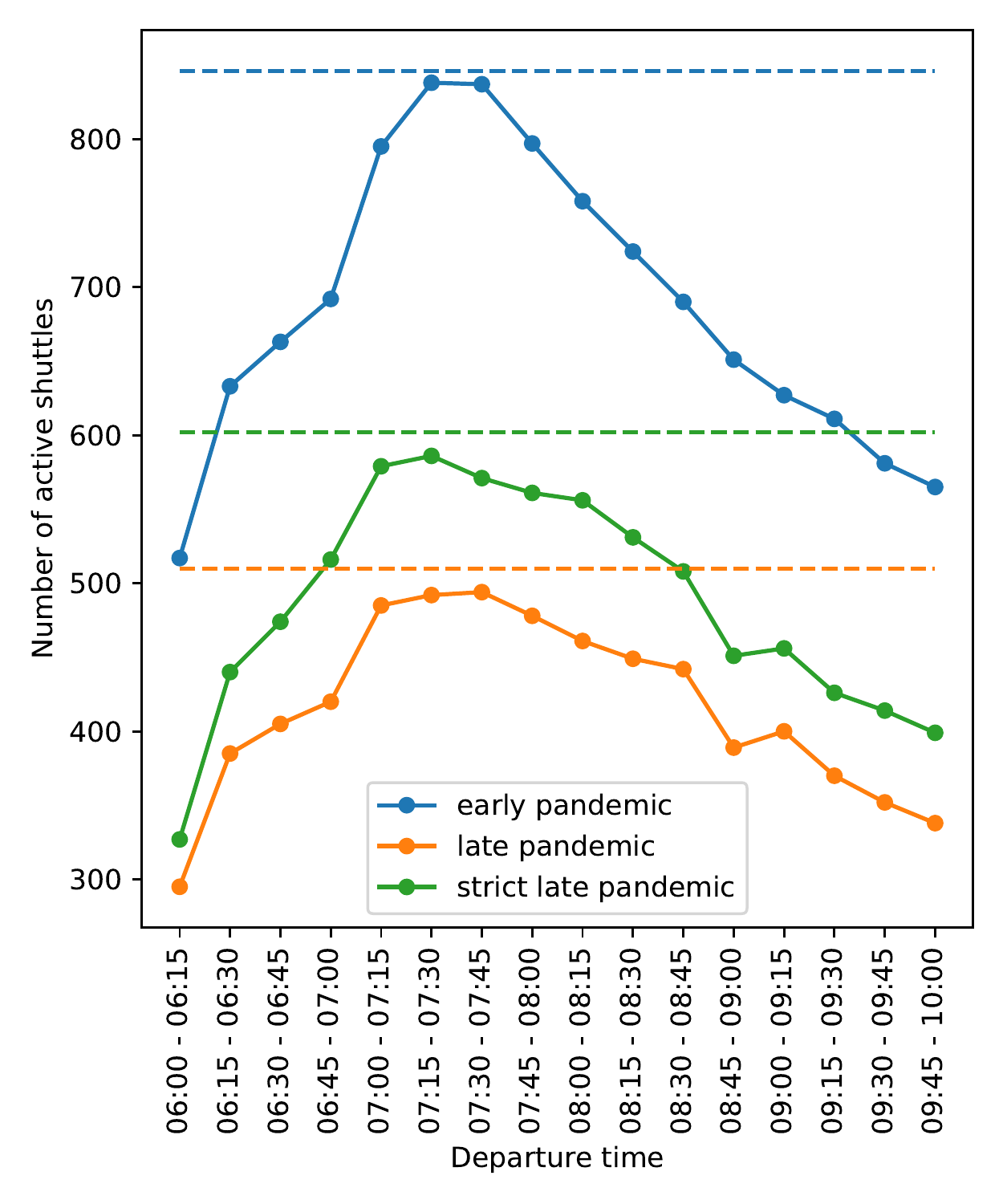}
    \captionof{figure}{Number of Shuttles Active over Time for Different Pandemic Scenarios (Number Available as Dashed Line).}
    \label{fig:multi:time:numberactive}
\end{figure}

In the baseline scenario, the capacity of the buses and the trains is
not restrictive, and this still holds true for the pandemic scenarios.
Table~\ref{tab:pandemic_overview} shows that the average bus
utilization does not exceed 17\%, while 50\% of capacity is available.
The average train utilization does not exceed 13\% in any pandemic
scenario.  This allows for additional social distancing, and it
provides opportunities for taking some vehicles out of the rotation
and cleaning them during the day, without creating capacity shortages.
A simulation of the early pandemic scenario confirms that halving the rail frequency would double the rail waiting time \changed{and the train utilization} as expected, but does not otherwise affect the system.
\changed{This also suggests that there is an opportunity to save costs during a pandemic by modifying the rail system.
This is studied for traditional systems by \citet{GkiotsalitisCats2021-OptimalFrequencySetting}, for example, and is also an interesting direction for future research in the context of ODMTS.}
Figure~\ref{fig:multi:time:numberactive} evaluates the utilization of
the shuttles during the pandemic.  For each pandemic scenario, it
shows the number of active shuttles during the day, compared to the
available number.  For all scenarios, the number of shuttles is
sufficient, and most of the shuttles are being used during the peak
from 7am to 8am.

While the performance of the ODMTS is similar for all scenarios, it
should be verified that the system remains affordable.  The budget and
cost columns in Table~\ref{tab:pandemic_overview} compare the cost of
the system to the reduced budget that results from decreased fare
revenues.
The cost is reported both with and without additional cleaning cost for sanitizing and disinfecting the fleet during the COVID-19 pandemic.
The cleaning cost is estimated at \$3.37 per hour for buses, and \$1.69 per hour for shuttles (see Appendix~\ref{app:costs}).
For the early pandemic scenario, the cost of the system (\$99k excl., \$105k incl.) stays within the \$106k budget, even when COVID-19 fleet sanitation is included.  This suggests that operating
the ODMTS at the same service level during this stage of the pandemic
is economically viable.  Section~\ref{sec:case_study:largerreductions}
explores how the quality of service is affected when a smaller budget
is enforced.  For the late and strict late pandemic scenarios, the
budget far exceeds the cost.

\changed{The resiliency of ODMTS for a fixed network design stems from the flexibility of the on-demand shuttles: Scaling the number of shuttles appropriately can compensate for large changes in demand and reductions in capacity, without exceeding the budget.}
ODMTS differs in this respect from
traditional systems, for which providing the same level of service
within a reduced budget is often impossible.

\subsection{Resiliency under Redesign}
\label{sec:case_study:redesign}

\begin{figure}[!t]
    \centering
    \includegraphics[trim=40cm 1cm 40cm 0cm,clip, height=0.4\textheight]{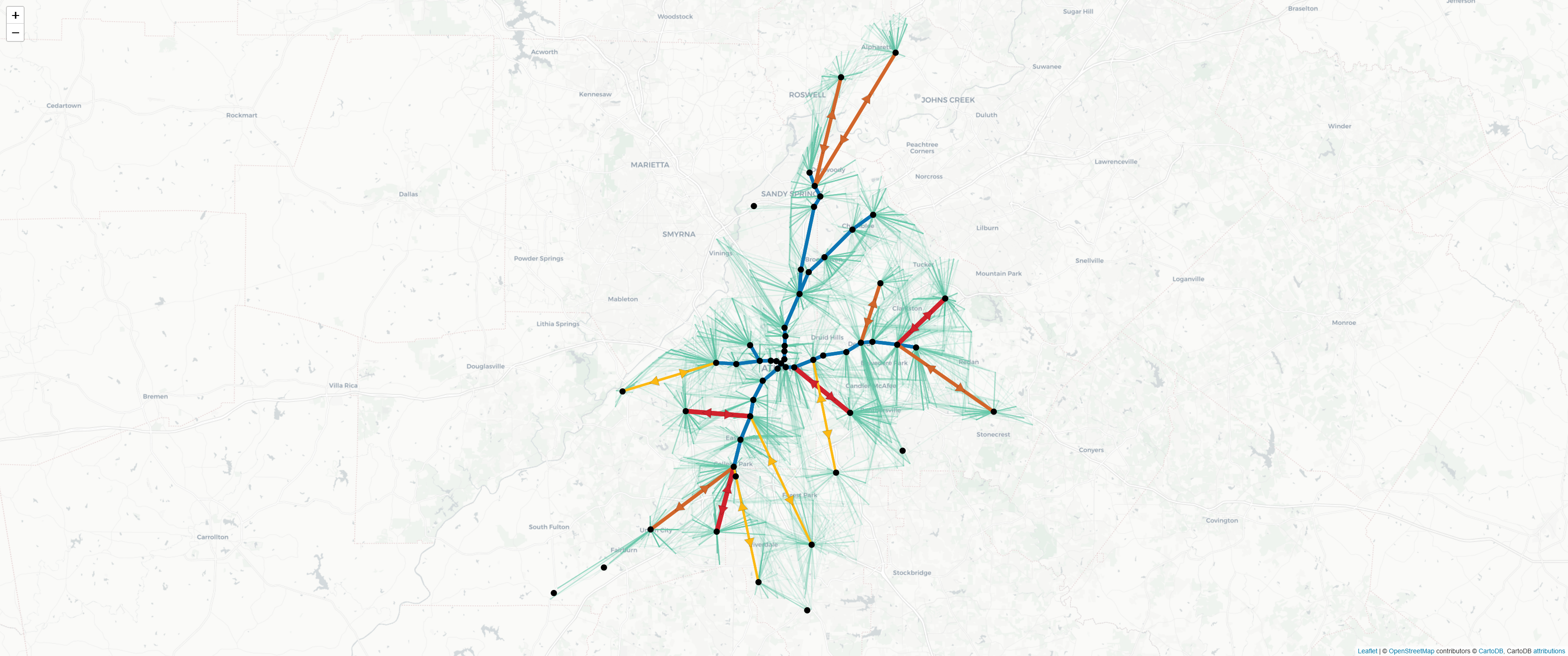}
    \caption{ODMTS Redesign for Early Pandemic Ridership.}
    \label{fig:early_pandemic_redesign:design}
\end{figure}
\begin{figure}[!t]
    \centering
    \includegraphics[trim=40cm 1cm 40cm 0cm,clip, height=0.4\textheight]{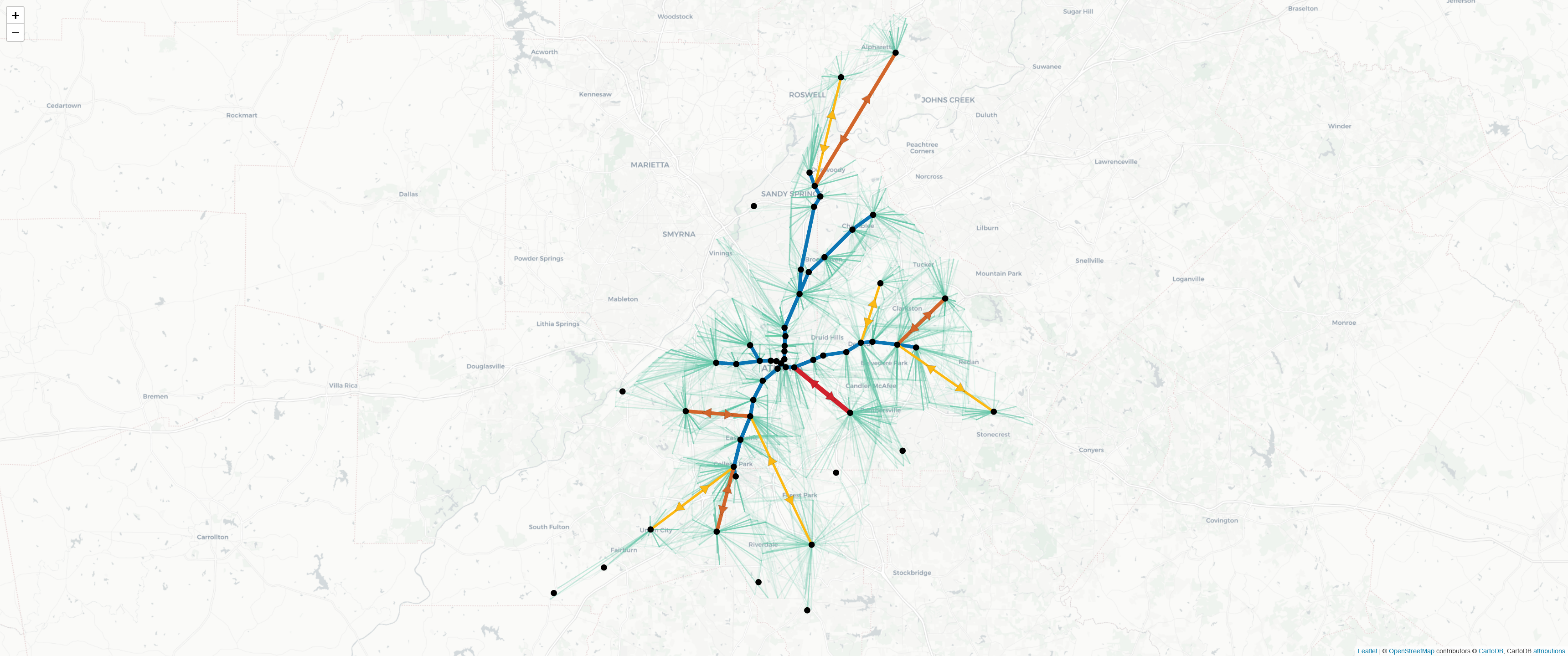}
    \caption{ODMTS Redesign for Late Pandemic Ridership.}
    \label{fig:late_pandemic_redesign:design}
\end{figure}

This section considers the cases where the transit agency redesigns the
network to better fit the pandemic circumstances and obtain additional
savings, especially if the transit system is expected to be affected for a longer
period of time.
For each pandemic scenario, the ODMTS design problem is solved for the respective pandemic ridership level.
To account for cleaning and social distancing, the bus cost per hour $c^{bus}$ is increased by \$3.37, and the shuttle cost per mile $\bar{c}^{shuttle}$ is replaced by the actual cost (incl. cleaning) obtained from the simulations under fixed design.
Figures~\ref{fig:early_pandemic_redesign:design} and
\ref{fig:late_pandemic_redesign:design} present the networks for the early and late
pandemic ridership levels respectively.

As ridership decreases, some
of the bus lines are closed, and others are reduced in frequency.  The
number of buses is reduced from 24 in the baseline design to 21 in the
early pandemic scenario and 13 for the late pandemic scenario.
Due to the relatively high cost of shuttles during a pandemic, the redesigns maintain relatively large bus subsystems compared to their demand.
Passengers are assigned \changed{fewer} direct shuttle trips, which allows the shuttle fleet to stay roughly the same size, despite the smaller number of buses.
Compared to 
the fixed design, the redesigned network in the early pandemic
scenario decreases the number of direct shuttle trips by 7.8\% and increases peak average bus utilization from 17\% to 20\%.
For the late pandemic scenario, direct shuttle trips decrease by 6.3\%, while bus utilization increases from 10\% to 18\%.

Redesigning the network results in \$2.0k (1.9\%) savings for the early
pandemic scenario, and \$2.9k (4.3\%) savings for the late pandemic
scenario.  In terms of performance, the average trip duration is
increased by less than a minute.  This is possible because the
the smaller networks fit better with the reduced demand, and passengers make less use of shuttles, which are more expensive during the pandemic due to cleaning cost and social distancing.

Redesigning the network based on the scenario may lead to additional
cost savings, without compromising performance.  This comes at the
cost of a more intensive use of buses, and slightly longer waiting
times for passengers using the bus subsystem.  The cost savings are
more significant in the late pandemic scenario, where ridership is
more depressed.
The increase in average bus wait time is under 2.5 minutes.

\changed{The possibility to redesign the network gives the ODMTS additional flexibility to deal with unexpected changes in demand, on top of the flexibility provided by the on-demand shuttles.
Rather than changing the design over time, the network design model may also be generalized to create a single fixed design with good performance in expectation.
This is an interesting direction for future research, but would require predicting the path of the pandemic, which has been very difficult for COVID-19.}

\subsection{Larger Budget Reductions}
\label{sec:case_study:largerreductions}

Section~\ref{sec:case_study:fixed} showed that, in the early pandemic
scenario, the cost of operating the ODMTS is within the reduced
budget.  This conclusion is critically dependent on how ridership
affects the budget.  In Atlanta, fares and direct revenues cover 33\%
of the operating expenses, but this number may be higher in other
areas.  In New York City, for example, the share is 55\%, and in the
San Francisco Bay Area it is even 70\%
\citep{FTA2018-NationalTransitDatabase}.  In these cases, loss of
ridership has a larger impact on the budget, and operating with the
proposed number of shuttles is no longer feasible.

\begin{figure}[!t]
    \centering
    \includegraphics[width=0.8\linewidth]{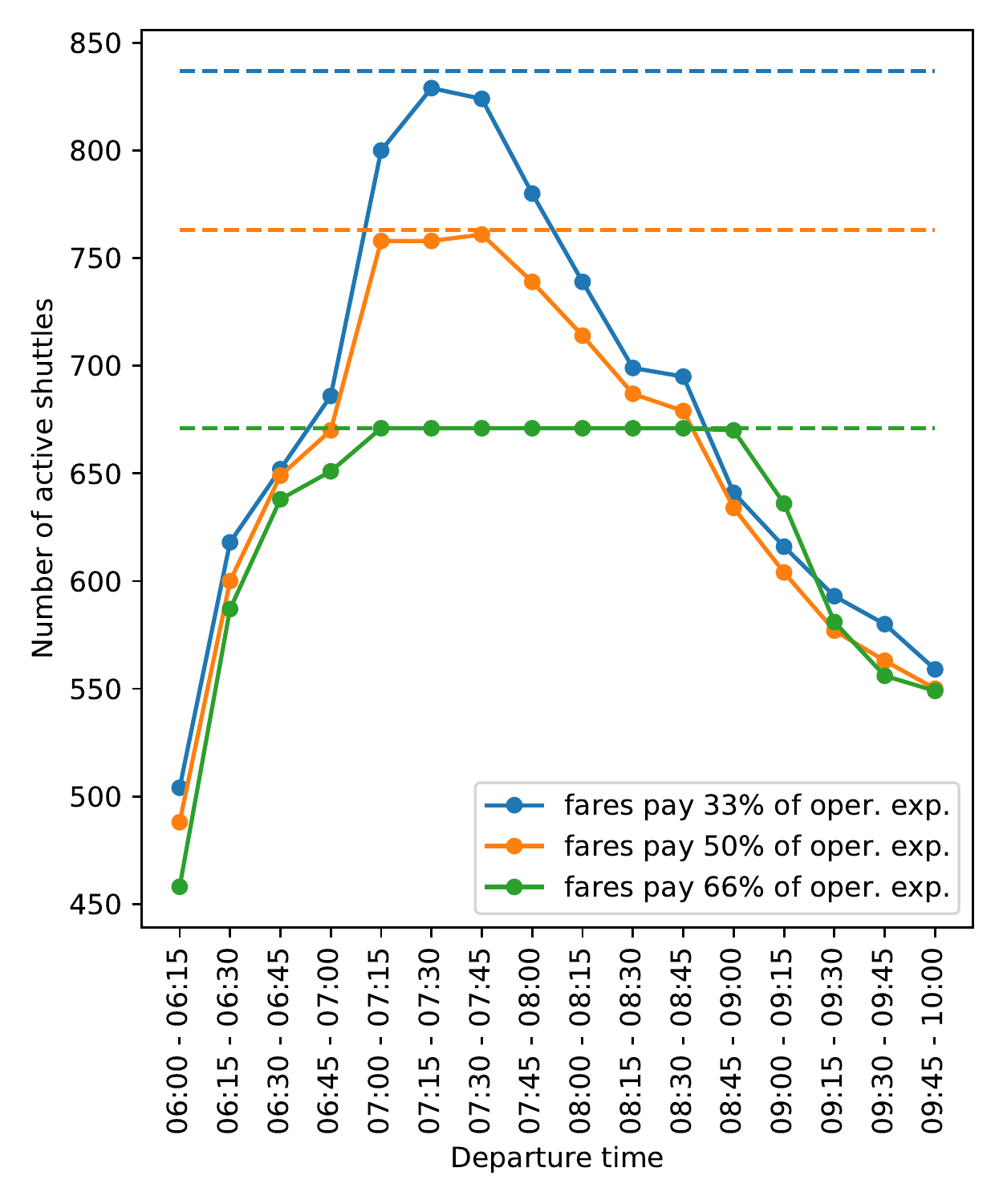}
    \captionof{figure}{Number of Shuttles Active over Time for the Early Pandemic Scenario for Different Budgets (Number Available as Dashed Line).}
    \label{fig:multi:time:numberactive:budget}
\end{figure}

To study the effect of larger budget reductions, two variations of the
early pandemic scenario are considered: one where fare revenues cover
50\% of the operating expenses, which corresponds to doubling the
current ticket price, and one where the share is doubled to 66\%. As
a result, the budget is reduced from \$106k to \$95k and \$84k
respectively.
Even using the redesigned network, the reduced budgets do not allow for operating 837 shuttles, as determined by the ridesharing and fleet-sizing
algorithms.  Instead, the number of shuttles is reduced to 763 and 671
respectively, to fit the
budget, including additional costs for cleaning. Figure~\ref{fig:multi:time:numberactive:budget} shows how
shuttle utilization is impacted by reducing the number of shuttles.
In the 33\% case, most shuttles are used during
the 7am to 8am peak, but there are always shuttles available for
incoming requests.  When the number of shuttles is reduced to 763,
almost all shuttles are active at the same time for about an hour.
The system does not seem to be overwhelmed, however, as the number of
active shuttles decreases quickly after the peak.  The 66\% case
paints a different picture: to operate the system with only 671
shuttles, the fleet is completely occupied for about two hours to
keep up with incoming requests.

\begin{figure}[!t]
    \centering
    \includegraphics[width=0.8\linewidth]{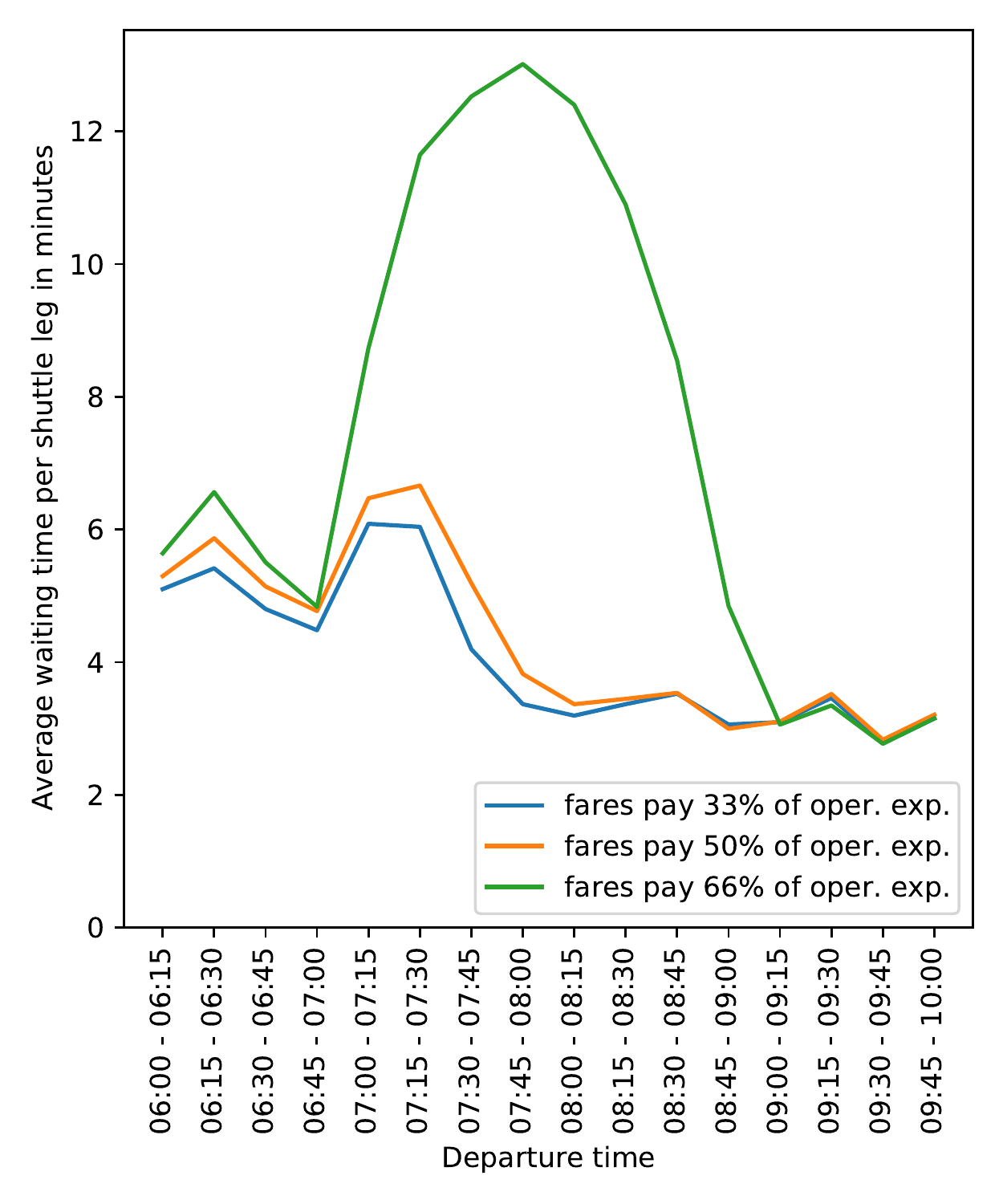}
    \captionof{figure}{Average Shuttle Waiting Time over Time for Different Budgets.}
    \label{fig:multi:time:waitingshuttle}
\end{figure}

When the fleet is busy, this affects the shuttle waiting time.
Figure~\ref{fig:multi:time:waitingshuttle} compares the waiting times
for a shuttle leg under the different budgets.  It can be seen that
reducing the budget results in a higher peak in shuttle waiting times.
Furthermore, the peak is wider at lower budgets, corresponding to the busy periods in
Figure~\ref{fig:multi:time:numberactive:budget}.  Reducing the number
of shuttles from 837 (33\% case) to 763 (50\% case) only has a minor
impact on performance.  The waiting time for a shuttle leg increases
by less than one minute during the peak, and the average trip
duration for all passengers remains practically the same.  Going down
to 671 shuttles (66\% case) presents a more significant trade-off
between cost and convenience.  The impact on the average trip duration
is only two minutes, but the average shuttle leg waiting time goes up
to 13 minutes during the peak.  It is up to the transit agency to
decide whether this wait is acceptable during a pandemic, given that
it allows for operating the system within a very restrictive budget that also covers additional cleaning,
without compromising access to public transit.

\section{\changed{Conclusions}}
\label{sec:conclusion}

This paper studied the resiliency of ODMTS during a pandemic, which
\changed{comes with} depressed ridership and revenues, and increased safety
requirements.  Operating in this environment is a challenge for public
transit agencies, which are forced to cut routes, reduce service
hours, or lower the service frequency.
An ODMTS can respond to a pandemic by scaling down the number of shuttles as needed. However, it was an open issue whether, under such \changed{conditions}, the ODMTS can still provide accessible public transit with a high quality of service that stays within the more restricted budget.
To answer this \changed{question}, the paper presented an ODMTS optimization pipeline which was applied to a variety of pandemic scenarios and evaluated using a high-fidelity simulator.

The case study considered the transit system for the city of Atlanta
based on actual data and observations from the COVID-19 pandemic.
The starting point was a baseline ODMTS designed under regular, non-pandemic circumstances. During a pandemic, the rail and
bus schedules may be kept the same, while the number of shuttles is
scaled down.  It was shown that this allows the ODMTS to provide
almost identical performance during a pandemic.  Furthermore, the cost
of the system decreases sufficiently to make up for lost revenues.
Transit agencies may also decide to redesign the network during a
pandemic to better fit the reduced demand.  The results showed that
this allows for additional cost savings without compromising service
quality, and this benefit is more significant when demand is more
depressed.  Finally, operating the ODMTS under a stricter budget was
considered to account for the realities of other transit systems.  The
simulations showed a trade-off between the budget and the peak shuttle
waiting times, but overall the average waiting time is not much
affected.

The case study results suggest that ODMTS provide significant
resiliency during pandemic response.  Furthermore, it supports earlier
findings that combining on-demand shuttles with high-frequency buses
and trains may provide a cost-efficient alternative to traditional
public transit systems.


\subsection*{Acknowledgments}
This research is partly supported by NSF Leap HI proposal NSF-1854684
and Department of Energy Research Award 7F-30154.

\bibliographystyle{trc}
\bibliography{references}


\appendix
\section{Cost Calculations}
\label{app:costs}
To compare the cost of the current system to that of an ODMTS, it is necessary to estimate the cost of operating a bus compared to operating an on-demand shuttle.
MARTA reports that the operating expenses for buses are \$104.10 per vehicle revenue hour.
However, this number includes expenses for general administration and facility maintenance, and excludes depreciation of the vehicles.
To present a fair comparison, it is assumed that buses and shuttles only differ in labor cost, vehicle maintenance cost, and vehicle depreciation.
In particular, the cost of fuel and lubricants is ignored, as this is a relatively small component of the cost.
If it were to be included, this would favor shuttles due to their better fuel efficiency.

The cost of one person-hour of operating a bus is \$23.52 in salaries and wages, to which \$16.03 (68\%) is added in fringe benefits \citep{FTA2018-NationalTransitDatabase}.
It is assumed that one vehicle revenue hour requires one person-hour of labor, which results in a conservative labor cost estimate.
The maintenance cost is estimated at \$19.17 per vehicle revenue hour.
This number results from adding up salaries and wages, fringe benefits, and materials in the `vehicle maintenance' category, and dividing by the number of vehicle revenue hours.
The majority of the MARTA fleet consists of 40-foot compressed natural gas buses with a current purchase price of \$625k \citep{FTA2018-NationalTransitDatabase,Dickens2020-PublicTransportationVehicle}.
These buses have a useful life benchmark of 12 years, and are used for 3878 vehicle revenue hours per year on average.
This corresponds to a \$13.43 depreciation per revenue hour.

For the purpose of comparing buses to shuttles, this brings the cost of a bus to \$72.15 per hour.
MARTA operates 465 buses in maximum service \citep{FTA2018-NationalTransitDatabase}, which corresponds to \$134k during the morning peak.

It is assumed that shuttle drivers are paid the average wage of a bus driver in Atlanta and that they receive the same 68\% of fringe benefits.
This amounts to \$14.42 per hour in wages \citep{Glassdoor2020-BusDriverSalaries} and \$9.83 in fringe benefits, for a total labor cost of \$24.25 per hour.
The maintenance cost is estimated at \$0.09 per mile for a medium sedan \citep{Association2019-YourDrivingCosts}, which results in \$1.13 per hour, assuming that shuttles and buses are used for the same number of hours per day and have the same daily mileage on average.
For the depreciation cost, it is assumed that shuttles have a purchase price of \$30k and are replaced every four years, which is after 191k miles.
This leads to a \$1.93 depreciation cost per hour.
Adding up these numbers, the comparable total cost of a shuttle is \$27.31 per vehicle revenue hour.

From an economic perspective, the public transit agencies may aim to recruit ridesharing drivers instead of bus drivers for operating the shuttles.
According to the Ridester's Driver Earnings Survey, the average Uber driver in the US earns \$19.36 per hour \citep{Helling2020-Ridesters2020Independent}.
The car is the driver's property, and the driver is responsible for the associated costs.
Based on the maintenance and depreciation costs as above, this results in an effective hourly wage of under \$16.30, which is low compared to \$24.25 for a bus driver.

The additional cost for sanitizing and disinfecting vehicles during COVID-19 is estimated from data provided by Gwinnett County Transit, the public transit system in the neighboring county of Gwinnett.
The data spans March to November 2020.
The average cost per operating hour is calculated by dividing the total cost of fleet sanitation by the number of bus operating hours, which results in a cost of \$3.37 per hour.
This paper assumes that cleaning significantly smaller on-demand shuttles is possible for half that cost, that is, \$1.69 per hour.

\vfill

\end{document}